\documentclass[12pt]{article}
\usepackage{textcomp}
\usepackage{graphicx}           
\usepackage{epstopdf}
\usepackage{color}              
\usepackage{xcolor,graphicx}
\usepackage{geometry}
\usepackage{float}
\usepackage{times}
\usepackage{indentfirst}        
\usepackage{amsmath,amssymb,bm} 
\usepackage{amsthm}
\usepackage{cases}              
\usepackage{setspace}
\usepackage{hyperref}
\usepackage{titlesec}
\usepackage[all]{xy}            
\usepackage{tikz}
\usepackage{multirow}
\usepackage{graphicx} 
\usepackage{epstopdf}
\usepackage{caption}
\usepackage{diagbox}
\usepackage{mathrsfs}
\usepackage{tikz}
\usepackage{dsfont}
\usepackage{graphics}
\usepackage{hyperref,cleveref,cite}
\usepackage{tcolorbox}
\usepackage{fancyhdr,fancyvrb}
\usepackage{xcolor,graphicx,geometry,float,indentfirst,cases,setspace,titlesec}
\usepackage{diagbox}

\newtheorem{thm}{Theorem}[section]
\newtheorem{defi}[thm]{Definition}
\newtheorem{lem}[thm]{Lemma}
\newtheorem{prop}[thm]{Proposition}
\newtheorem{coro}[thm]{Corollary}

\newtheorem{pf}{proof}[section]
\theoremstyle{remark}

\begingroup
\makeatletter \@addtoreset{equation}{section} \makeatother
\makeindex 
\def\qed{\hfill \rule{4pt}{7pt}}
\def\pf{\vskip 0.2cm {\noindent \bf Proof.}\quad}

\begin{document}

\begin{center}

 {\Large \bf Some Refinements of Stanley's Shuffle Theorem}

\end{center}

\begin{center}
{Kathy Q. Ji}$^{1}$ and {Dax T.X. Zhang}$^{2}$  \vskip 2mm

$^{1,2}$ Center for Applied Mathematics,  Tianjin University, Tianjin 300072, P.R. China\\[6pt]
   \vskip 2mm

    $^1$kathyji@tju.edu.cn and $^2$zhangtianxing6@tju.edu.cn
\end{center}

\vskip 6mm \noindent {\bf Abstract.}  We give a combinatorial proof of Stanley's shuffle theorem by using the insertion lemma of Haglund, Loehr and Remmel. Based on this combinatorial construction, we establish several refinements of  Stanley's shuffle theorem.

\noindent
{\bf Keywords:} descent, major index, permutation, shuffle, partition

\noindent
{\bf AMS Classification:} 05A05, 05A19, 11P81

 \vskip 6mm

\section{Introduction}

Stanley's shuffle theorem gives an explicit expression for the generating function of the number of shufflings of two disjoint permutations $\sigma$ and $\pi$ with  $k$ descents and the major index being $t$.  Let us recall some common notation and terminology on permutations as used in \cite[Chapter 1]{Stanley-2012}. We say that $\pi=\pi_1\cdots \pi_n$ is a permutation of length $n$   if it is a sequence of $n$ distinct letters--not necessarily from $1$ to $n$. Let  $\ell(\pi)$ denote the length of the permutation $\pi$. For example, $\pi=9\,3\,8\,10\,12\,3\,7$ is a permutation of length $7$, and so $\ell(\pi)=7$. Let $\mathfrak{S}_n$ denote the set of all permutations of length $n$.  We say that $1\leq i \leq n-1$ is a {\rm des}cent of $\pi \in \mathfrak{S}_n$ if $\pi_i>\pi_{i+1}$ and $1\leq i \leq n-1$ is an ascent of $\pi \in \mathfrak{S}_n$ if $\pi_i<\pi_{i+1}$.  The set of  descents of $\pi$ is called the descent set of $\pi$, denoted ${\rm Des}(\pi)$
and  the number of its descents is called the descent number, denoted ${\rm des}(\pi)$. The   major index of $\pi$, denoted ${\rm maj}(\pi)$,  is defined to be the sum of its  descents. To wit,
\[{\rm maj}(\pi):=\sum_{k \in {\rm Des}(\pi)}k.\]
Let $d_k(\pi)$ denote the number of  descents in $\pi$ greater than or equal to $k$. Then,
\[d_1(\pi)={\rm des}(\pi) \quad \text{and} \quad  {\rm maj}(\pi)=\sum_{k=1}^n d_k(\pi).\]

Let $\sigma \in \mathfrak{S}_n$ and $\pi \in \mathfrak{S}_m$ be two disjoint permutations, that is, permutations with no letters in common. We say that $\alpha \in \mathfrak{S}_{n+m}$ is a shuffle of $\sigma$ and $\pi$ if both $\sigma$ and $\pi$  are
subsequences of $\alpha$. The set of shuffles of  $\sigma$ and $\pi$  is denoted $\mathfrak{S}(\sigma, \pi)$. For example,
\begin{align*}
    \mathfrak{S}(2\,6\,3,1\,4)=\{2&\,6\,3\,1\,4, 2\,6\,1\,3\,4, 2\,6\,1\,4\,3, 2\,1\,4\,6\,3, 2\,1\,6\,3\,4, \,
2\,1\,6\,4\,3, 1\,2\,4\,6\,3,  \\[2pt] &1\,4\,2\,6\,3,\, 1\,2\,6\,3\,4, 1\,2\,6\,4\,3\}.
\end{align*}
It is easy to see that the number of permutations in $\mathfrak{S}(\sigma, \pi) $ is ${\ell(\sigma)+\ell(\pi) \choose \ell(\pi)}.$

Define
\[ {S}_{k,q}(\sigma,\pi)=\sum_{\alpha\in \mathfrak{S}(\sigma,\pi) \atop {\rm des}(\alpha)=k}q^{{\rm maj}(\alpha)}.\]
In light of the $q$-Pfaff-Saalsch\"utz identity in his setting of $P$-partitions,  Stanley \cite{Stanley-1972}  obtained a compact expression for ${S}_{k,q}(\sigma,\pi)$ in terms of the Gaussian polynomial (also  called the $q$-binomial coefficients), as given by
\[{n+m \brack n}=\begin{cases}
\small \displaystyle\frac{(1-q^{n+m})(1-q^{n+m-1})\cdots (1-q^{m+1})}{(1-q^n)(1-q^{n-1})\cdots (1-q)}, & \text{for } n, m\geq 0,\\[5pt]
0,& \text{otherwise},
\end{cases}\]
for $n,m$ are non-negative integers, see Andrews \cite[Chapter 1]{Andrews-1976}.

\begin{thm}[Stanley's shuffle theorem]\label{stanley}
Let $\sigma$ and $\pi$ be two disjoint permutations.  Then
\begin{align}
   \sum_{\alpha\in \mathfrak{S}(\sigma,\pi) \atop {\rm des}(\alpha)=k}q^{{\rm maj}(\alpha)}&= {\ell(\sigma)-{\rm des}(\sigma)+{\rm des}(\pi) \brack k-{\rm des}(\sigma)} {\ell(\pi)-{\rm des}(\pi)+{\rm des}(\sigma)  \brack  k-{\rm des}(\pi)} \nonumber \\[5pt]
   &\quad \quad \quad \quad\times q^{{\rm maj}(\sigma)+{\rm maj}(\pi)+(k-{\rm des}(\pi))(k-{\rm des}(\sigma))}.
\end{align}
\end{thm}

Stanley asked for a proof of Theorem \ref{stanley} which avoids the use of the $q$-Pfaff-Saalsch\"utz identity (see \cite[Eq.3.3.11]{Andrews-1976}). Bijective proofs of  Stanley's shuffle theorem have been given by Goulden \cite{Goulden-1985}  and Stadler \cite{Stadler-1999}. Goulden's proof is obtained by finding bijections for lattice path representations of shuffles which reduce  $\sigma$ and $\pi$ to canonical permutations, for which  the generating function is easily given. Stadler's bijection is more elementary, but the inverse of Stadler's map is not very explicit. In this paper, we first give an explicit bijective proof of Theorem \ref{stanley} by using the insertion lemma of Haglund, Loehr and Remmel \cite{Haglund-Loehr-Remmel-2005}. It turns out that the insertion lemma of Haglund, Loehr and Remmel is equivalent to Stanley's shuffle theorem in the case   $\ell(\pi)=1$. It should be mentioned that Novick \cite{Novick-2010} used the insertion lemma of Haglund, Loehr and Remmel to give a bijective proof of the following theorem due to Garsia and Gessel \cite{Garsia-Gessel-1979}.

\begin{thm}[Garsia and Gessel]\label{Garsia-Gessel}
Let $\sigma$ and $\pi$ be two disjoint permutations.   Then
\begin{equation}
    \sum_{\alpha \in \mathfrak{S}(\sigma,\pi)}q^{{\rm maj}(\alpha)}=
   {\ell(\sigma)+\ell(\pi) \brack \ell(\pi)}q^{{\rm maj}(\sigma)+{\rm maj}(\pi)} .
\end{equation}
\end{thm}

In fact, Theorem \ref{Garsia-Gessel} can be derived from Theorem \ref{stanley} by employing $q$-analogue of the Chu-Vandermonde summation (see \cite[Eq.3.3.10]{Andrews-1976}),
\[\sum_{k=0}^h {n \brack k}{m \brack h-k}q^{(n-k)(h-k)}={m+n \brack h}.\]

Based on this combinatorial construction, we obtain four refinements of Stanley's shuffle theorem, see Theorem \ref{stanley-r2}, Theorem \ref{stanley-r1}, Theorem \ref{stanley-r3}  and Theorem \ref{stanley-r4}  in Section 4. As immediate consequences of these four refinements, we obtain two more general refinements.
More precisely,  let $\mathfrak{S}^{la,\, sa,\,sb}(\sigma,\pi)$ denote the set of shuffles $\alpha=\alpha_1\cdots\alpha_{\ell(\alpha)}$ of two disjoint permutations $\sigma=\sigma_1\cdots \sigma_{\ell(\sigma)}$ and $\pi=\pi_1\cdots\pi_{\ell(\pi)}$ such that $\alpha_{\ell(\alpha)}=\min\{\sigma_{\ell(\sigma)}, \pi_{\ell(\pi)}\}$, satisfying the conditions in Definition \ref{defi-refi-1} and Definition \ref{defi-refi-3}. For example,
\[\mathfrak{S}^{la,\,sa,\,sb}(2\,6\,3,1\,4)=\{ 2\,1\,6\,4\,3\}.\]

Combining Theorem \ref{stanley-r2}, Theorem \ref{stanley-r1} and Theorem \ref{stanley-r4},
we arrive at the following refinement of Stanley's shuffle theorem.

\begin{thm}
Let $\sigma$ and $\pi$ be two disjoint permutations. Then
\begin{align*}
     \sum_{\alpha\in \mathfrak{S}^{la,\, sa,\,sb}(\sigma,\pi)\atop  {\rm des}(\alpha)=k}q^{{\rm maj}(\alpha)}& = {\ell(\sigma)-{\rm des}(\sigma)+{\rm des}(\pi)-2 \brack k-{\rm des}(\sigma)-2}{\ell(\pi)-{\rm des}(\pi)+{\rm des}(\sigma)-1 \brack k-{\rm des}(\pi)-1} \\[5pt]
    &\quad \quad \times q^{{\rm maj}(\sigma)+{\rm maj}(\pi)+\ell(\sigma)+\ell(\pi)-1+(k-{\rm des}(\pi)-1)(k-{\rm des}(\sigma)-1)}.
\end{align*}
\end{thm}

Let $\mathfrak{S}^{la,\, lb,\,sb}(\sigma,\pi)$ denote the set of shuffles $\alpha=\alpha_1\cdots\alpha_{\ell(\alpha)}$ of two disjoint permutations $\sigma=\sigma_1\cdots \sigma_{\ell(\sigma)}$ and $\pi=\pi_1\cdots\pi_{\ell(\pi)}$ such that $\alpha_{\ell(\alpha)}=\min\{\sigma_{\ell(\sigma)}, \pi_{\ell(\pi)}\}$, satisfying the conditions in Definition \ref{defi-refi-1} and Definition \ref{defi-refi-4}. For example,
\[\mathfrak{S}^{la,\, lb,\,sb}(2\,6\,3,1\,4)=\{ 1\,2\,6\,4\,3\}.\]

Using Theorem \ref{stanley-r2} and Theorem \ref{stanley-r1}, as well as Theorem \ref{stanley-r3},
we obtain the following refinement of Stanley's shuffle theorem.

\begin{thm}
Let $\sigma$ and $\pi$ be two disjoint permutations. Then
\begin{align*}
     \sum_{\alpha\in \mathfrak{S}^{la,\, lb,\,sb}(\sigma,\pi)\atop  {\rm des}(\alpha)=k}q^{{\rm maj}(\alpha)}& = {\ell(\sigma)-{\rm des}(\sigma)+{\rm des}(\pi)-1 \brack k-{\rm des}(\sigma)-1}{\ell(\pi)-{\rm des}(\pi)+{\rm des}(\sigma)-2 \brack k-{\rm des}(\pi)-1} \\[5pt]
    &\quad \quad \times q^{{\rm maj}(\sigma)+{\rm maj}(\pi)+\ell(\sigma)+\ell(\pi)-1+(k-{\rm des}(\pi)-1)(k-{\rm des}(\sigma))}.
\end{align*}
\end{thm}

To conclude the introduction, let us say a few words on the recent work  that has built upon Stanley's shuffle theorem.
Inspired by Stanley's shuffle theorem, Gessel and Zhuang \cite{Gessel-Zhuang-2018} introduced the concept of
shuffle compatible and stated that Stanley's shuffle
theorem  imply that maj and des are shuffle compatible.
Gessel and Zhuang \cite{Gessel-Zhuang-2018} further
investigated the shuffle compatibility property of
other permutation statistics  and posed several
conjectures involving   the shuffle compatibility. Some
of these conjectures were confirmed by    Baker-Jarvis
and Sagan \cite{Baker-Jarvis-Sagan-2020}, Grinberg
\cite{Grinberg-2018} and   Yang and Yan \cite{Yang-Yan-2022}.  In particular, Baker-Jarvis and
Sagan \cite{Baker-Jarvis-Sagan-2020} provided  unified
bijective techniques to give a demonstration of shuffle
compatibility. Cyclic shuffle and cyclic shuffle
compatibility were  investigated by Adin, Gessel, Reiner and Roichman \cite{Adin-Gessel-Reiner-Roichman-2021} and Domagalski, Liang, Minnich and Sagan \cite{Domagalski-Liang-Minnich-Sagan-2020}, respectively. In
\cite{{Ji-Zhang-2022}},  we established a cyclic analogue of Theorem \ref{stanley}.

Last but not least, we would like to mention one interesting consequence of Stanley's shuffle theorem, which was   explicitly stated by Sagan and Savage \cite[Corollary 2.4]{Sagan-Savage-2012} and proved by using Foata's fundamental bijection \cite{Foata-1968}.
\begin{thm}[Sagan and Savage]\label{ss} For $m,n\geq 1$,
\begin{align}
   \sum_{\alpha\in \mathfrak{S}(1^m,\,2^n) \atop {\rm des}(\alpha)=k}q^{{\rm maj}(\alpha)}&= {m \brack k} {n  \brack  k}   q^{k^2}.
\end{align}
\end{thm}
It should be noted that a specialization of our combinatorial construction provides an alternative proof of  Theorem \ref{ss}.

The paper is organized as follows. In Section 2, we recall the insertion lemma of Haglund, Loehr and Remmel and demonstrate that the insertion lemma of Haglund, Loehr and Remmel is equivalent to Stanley's shuffle theorem in the case   $\ell(\pi)=1$. Section 3 is devoted to the bijective proof of Theorem \ref{stanley} based on the insertion lemma of Haglund, Loehr and Remmel. In Section 4, we establish four refinements of Stanley's shuffle theorem relying on this combinatorial construction.

\section{The insertion lemma of Haglund, Loehr and Remmel}

This section is devoted to illustrating the insertion lemma of Haglund, Loehr and Remmel \cite{Haglund-Loehr-Remmel-2005}.
We follow the terminology, notation and the example in Section 4 of their paper.

Assume that  $\sigma=\sigma_1\cdots \sigma_n \in \mathfrak{S}_n$ and $r \not \in \sigma$ (that is, there does not exist $1\leq j\leq n$ such that $\sigma_j=r$). For $0\leq i\leq n$,  let $\sigma^{(i)}(r)$ denote the permutation of length $n+1$ obtained by inserting $r$ into $\sigma$ before $\sigma_{i+1}$. Here we assume that $\sigma^{(n)}(r)$ denotes the permutation of length $n+1$ obtained by inserting $r$ into $\sigma$ after $\sigma_n$. The insertion lemma of Haglund, Loehr and Remmel  \cite{Haglund-Loehr-Remmel-2005} showed that   no matter what  the relative value of $r$ is with respect to the elements in $\sigma$,
\begin{equation}\label{inserlem}
 \sum_{i=0}^{n}q^{{\rm maj}(\sigma^{(i)}(r))}=(1+q+\cdots+q^{n}) q^{{\rm maj}(\sigma)}.
\end{equation}
This relation can be used to establish the following celebrated formula due to MacMahon \cite{MacMahon-1978}.
\begin{equation}\label{MacMahon}
  \sum_{\sigma \in \mathfrak{S}_n} q^{{\rm maj}(\sigma)}=[n]_q!,
\end{equation}
where
$[n]_q!=[1]_q[2]_q\cdots [n]_q$ with $[k]_q=1+q+\cdots+q^{k-1}$.

 Haglund, Loehr and Remmel \cite{Haglund-Loehr-Remmel-2005} also classified the possible spaces where they can insert $r$ into $\sigma$  into two sets called the right-to-left spaces which they denoted as RL-spaces and the left-to-right spaces which they
 denoted as LR-spaces. That is, a space $i$ is called a RL-space of $\sigma$ relative to $r$ if
 \begin{itemize}
     \item[1.] $i = n$ and $\sigma_n<r$,
 \item[2.] $i=0$ and $r<\sigma_1$,
 \item[3.] $0 < i < n$ and $\sigma_i > \sigma_{i+1}>r$,
\item[4.] $0 < i < n$ and $r>\sigma_i > \sigma_{i+1}$, or
\item[5.] $0 < i < n$ and $\sigma_i<r< \sigma_{i+1}$.
 \end{itemize}
A space $i$ is a LR-space of $\sigma$  relative to $r$ if it is not a RL-space of $\sigma$  relative to $r$. Assume that there are $l$ RL-spaces of $\sigma$ relative to $r$. Label the RL-spaces from
right to left with $0,\ldots, l-1$ and  label the LR-spaces from left to right with $l,\ldots, n$. We call this labeling the canonical labeling for $\sigma$ relative to $r$. For example, suppose that $r=5$ and $\sigma= 10 \ 1\  9\  8\  2\  7\  4\  3\  6$ is a permutation in $\mathfrak{S}_9$. By definition, we see the RL-spaces of $\sigma$ relative to $5$ are $0, 2, 3, 5, 7$ and $8$ and the LR-spaces of $\sigma$ relative to $5$
are $1, 4, 6 $ and $9$. The canonical labeling of $\sigma$ relative to $r$ is
\[_{\tiny  \bf 5}10 \, _{\tiny 6}1\,  _{\tiny  \bf 4}9\, _{\tiny  \bf 3}8\,  _{\tiny 7}2\, _{\tiny  \bf 2}7\,  _{\tiny 8}4\,  _{\tiny  \bf 1}3\, _{\tiny \bf 0}6\,  _{\tiny 9},\]
where the bold number in the subscript represents the labeling of the RL-spaces of $\sigma$ relative to $5$.

 Haglund, Loehr and Remmel \cite{Haglund-Loehr-Remmel-2005} established the following insertion lemma.

\begin{lem}[The insertion lemma]\label{hag-leo-rem}   Suppose that $\sigma=\sigma_1\cdots \sigma_n \in \mathfrak{S}_n$ and   $r\not \in \sigma$, and let $\sigma^{(i)}(r)$ denote the permutation obtained by inserting $r$ into $\sigma$ before $\sigma_{i+1}$. If the label at the $i$-th space in the cannonical labeling of $\sigma$ relative to $r$ is equal to $a$, then
\[{\rm maj}(\sigma^{(i)}(r) )=a+{\rm maj}(\sigma).\]
\end{lem}

With a careful examination of the definitions of the RL-spaces and the LR-spaces, we obtain  the following lemma, which is useful in the proof of Stanley's shuffle theorem.

\begin{lem} \label{hag-leo-rem-add} Let $\sigma$, $r$ and $\sigma^{(i)}(r)$  be given in Lemma \ref{hag-leo-rem}. If $i$ is a RL-space of $\sigma$ relative to $r$, then ${\rm des}(\sigma^{(i)}(r))={\rm des}(\sigma)$. If $i$ is a LR-space of $\sigma$ relative to $r$, then ${\rm des}(\sigma^{(i)}(r))={\rm des}(\sigma)+1$. Moreover, the number of   RL-spaces of $\sigma$ relative to $r$ is one more than the number of { des}cents in $\sigma$.
\end{lem}

\pf  Assume that there are $k$  descents in $\sigma$.  From the definitions of RL-spaces and LR-spaces, we find that if $i$ is a RL-space of $\sigma$ relative to $r$, then ${\rm des}(\sigma^{(i)}(r))={\rm des}(\sigma)=k$. Moreover, the  major increment  ${\rm maj}(\sigma^{(i)}(r) )-{\rm maj}(\sigma)=d_{i+1}(\sigma^{(i)}(r))$, and hence ${\rm maj}(\sigma^{(i)}(r) )-{\rm maj}(\sigma)\le k$. If $i$ is a LR-space of $\sigma$ relative to $r$, then ${\rm des}(\sigma^{(i)}(r))={\rm des}(\sigma)+1$, and it can be checked that the major increment  $k<{\rm maj}(\sigma^{(i)}(r) )-{\rm maj}(\sigma)\le n$.   Assume that there are $n_r$ RL-spaces of $\sigma$ and $n_l$  LR-spaces of $\sigma$. Clearly, $n_r+n_l=n+1$.
 By Lemma \ref{hag-leo-rem}, we see that the  major increment at each space of $\sigma$ is different, so we conclude that    $n_r\le k+1$ and $n_l\le n-k$. Since $n_r+n_l=n+1$, we derive that $n_r=k+1$ and $n_l=n-k$. This completes the proof.   \qed

We conclude this section with the proof of Stanley's shuffle theorem in the case $\ell(\pi)=1$ in view of Lemma \ref{hag-leo-rem} and Lemma \ref{hag-leo-rem-add}.

If  $\ell(\pi)=1$ in Theorem \ref{stanley} , then ${\rm des}(\pi)=0$  and ${\rm maj}(\pi)=0$, and
thus Theorem \ref{stanley}  reads as follows:
\begin{equation}\label{sst-sp}
   \sum_{\alpha\in \mathfrak{S}(\sigma,\pi) \atop {\rm des}(\alpha)=k}q^{{\rm maj}(\alpha)}= {\ell(\sigma)-{\rm des}(\sigma) \brack k-{\rm des}(\sigma)} {{\rm des}(\sigma)+1  \brack  k}  q^{{\rm maj}(\sigma)+k(k-{\rm des}(\sigma))}.
\end{equation}
Assume that $\ell(\sigma)=m$ and ${\rm des}(\sigma)=r$. In this case, we see that  the right-hand side of \eqref{sst-sp} is non-zero if and only if $k=r$ or $k=r+1$. Hence  \eqref{sst-sp} can be written as
\begin{equation}\label{sst-sp-a1}
   \sum_{\alpha\in \mathfrak{S}(\sigma,\pi) \atop {\rm des}(\alpha)=r}q^{{\rm maj}(\alpha)}=   {r+1  \brack r}  q^{{\rm maj}(\pi)}=(1+q+\cdots +q^r)\ q^{{\rm maj}(\sigma)}.
\end{equation}
and
\begin{equation}\label{sst-sp-a2}
   \sum_{\alpha\in \mathfrak{S}(\sigma,\pi) \atop {\rm des}(\alpha)=r+1}q^{{\rm maj}(\alpha)}=   {m-r  \brack 1}  q^{{\rm maj}(\pi)+r+1}=q^{r+1}(1+q+\cdots +q^{m-r-1})\ q^{{\rm maj}(\sigma)}.
\end{equation}
Obviously, the identity \eqref{sst-sp-a1} and the identity \eqref{sst-sp-a2} are immediate consequences of   Lemma \ref{hag-leo-rem} and Lemma \ref{hag-leo-rem-add}.

\section{The bijection}

In this section, we give a proof of Theorem \ref{stanley} in the general case with the aid of  Lemma  \ref{hag-leo-rem} and Lemma \ref{hag-leo-rem-add}. To state the proof, we need to recall   some  notation and terminology on partitions as  in \cite[Chapter 1]{Andrews-1976}.   {
A partition} $\lambda$ of a positive integer $n$ is a finite
nonincreasing sequence of nonnegative integers
$(\lambda_1,\,\ldots,\,\lambda_s)$ such that
$\sum_{i=1}^s\lambda_i=n.$  Then $\lambda_i$ are called the parts
of $\lambda$, where $\lambda_1$ is its largest part and $\lambda_s$ is its smallest part. The number of parts
of  $\lambda$ is called the length of $\lambda$, denoted
$\ell(\lambda).$  The weight of $\lambda$ is the sum of  parts of $\lambda$, denoted
 $|\lambda|.$   Let  $\mathcal{P}_{\leq n}(m)$  denote  the set of partitions $\lambda$ such that $\ell(\lambda)\leq n$ and $\lambda_1\leq m$. It is well-known that the Gaussian polynomial  has the following partition interpretation \cite[Theorem 3.1]{Andrews-1976}:
\begin{equation}\label{int-GassCoef}
{n+m \brack n}=\sum_{\lambda \in \mathcal{P}_{\leq n}(m)} q^{|\lambda|}.
\end{equation}
 In general, let  $\mathcal{P}_n(t,m)$  denote  the set of partitions $\lambda$ such that $\ell(\lambda)=n$, $\lambda_n\geq t$ and $\lambda_1\leq m$, we have
  \begin{equation}\label{int-GassCoeft}
q^{nt}{n+m-t \brack n}=\sum_{\lambda \in \mathcal{P}_n(t,m)} q^{|\lambda|}.
\end{equation}
When $t=0$, we see that   \eqref{int-GassCoeft} coincides  with  \eqref{int-GassCoef}.

Using   \eqref{int-GassCoeft}, we see that  Theorem \ref{stanley} is equivalent to the following combinatorial statement.

 \begin{thm}\label{Stanley-c} Suppose that $\sigma \in \mathfrak{S}_m$ and $\pi \in \mathfrak{S}_n$ are two disjoint permutations, where ${\rm des}(\sigma)=r$  and  ${\rm des}(\pi)=s$. Let $\mathfrak{S}(\sigma,\pi|k)$ denote the set of all shuffles of $\sigma$ and $\pi$ with $k$ {des}cents. Then there is a bijection $\Phi$ between $\mathfrak{S}(\sigma,\pi|k)$ and $ \mathcal{P}_{k-r}(k-s, m)\times \mathcal{P}_{n-k+r}(0,\,k-s)$, namely,   for   $\alpha \in \mathfrak{S}(\sigma,\pi|k)$, we have
 $(\lambda,\mu)=\Phi(\alpha)\in  \mathcal{P}_{k-r}(k-s, m) \times \mathcal{P}_{n-k+r}(0,\,k-s)$ such that
 \begin{equation} \label{thm-rel}
     {\rm maj}(\alpha)=|\lambda|+|\mu|+{\rm maj}(\sigma)+{\rm maj}(\pi).
 \end{equation}
 \end{thm}

 To prove Theorem \ref{Stanley-c}, we first give a {\rm des}cription of the map $\Phi$,
 and then we show that the map $\Phi$ is a bijection,  as {\rm des}ired in Theorem  \ref{Stanley-c}.

 \begin{defi}[The map $\Phi$] \label{defi-map} Let $\sigma=\sigma_1\cdots\sigma_m$ be a permutation   with $r$ {\rm des}cents and let $\pi=\pi_1\cdots\pi_n$ be a permutation  with $s$ {\rm des}cents.
 Assume that $\alpha=\alpha_1\cdots \alpha_{n+m}$ is  the shuffle of $\sigma$ and $\pi$ with $k$ {\rm des}cents.   The pair of partitions $(\lambda,\mu)=\Phi(\alpha)$ can be constructed as follows:  Let
$\alpha^{(i)}$ denote the permutation obtained by removing $\pi_1,\pi_{2},\ldots,\pi_{i}$ from $\alpha$. Obviously, $\alpha^{(n)}=\sigma$. Here we assume that  $\alpha^{(0)}=\alpha$.  For $1\leq i\leq n$, define
\begin{equation}\label{defi-t}
 t(i)={\rm maj}(\alpha^{(i-1)})-{\rm maj}(\alpha^{(i)})-d_i(\pi).
\end{equation}
 Since there are $k$ {\rm des}cents in $\alpha$ and there are $r$ {\rm des}cents in $\sigma$,  it follows that there exists $k-r$ permutations in $\alpha^{(1)},\ldots,\alpha^{(n)}$,  denoted $\alpha^{(i_1)},\ldots,\alpha^{(i_{k-r})}$ where $1\leq i_1<i_2<\cdots <i_{k-r}\leq n$,  such that ${\rm des}(\alpha^{(i_l-1)})={\rm des}(\alpha^{(i_l)})+1$  for $1\le l\le k-r$. Let $\{j_1,\ldots,j_{n-k+r}\} = \{1,\ldots,n\}
\backslash  \{i_1,i_2,\ldots, i_{k-r}\}$, where $1\leq j_1<j_2<\cdots<j_{n-k+r} \leq n$. Then  ${\rm des}(\alpha^{(j_l-1)})={\rm des}(\alpha^{(j_l)})$ for $1\le l\le n-k+r$.  The pair of partitions $(\lambda,\mu)=\Phi(\alpha)$ is defined by
\begin{equation}\label{defi-lambda}
   \lambda=(t(i_{k-r}),t(i_{k-r-1}),\ldots, t(i_{1})),
\end{equation}
and
\begin{equation}\label{defi-mu}
  \mu=(t(j_{1}),t(j_{2}),\ldots, t(j_{n-k+r})).
\end{equation}
\end{defi}

 For example, let \[\sigma=9\,3\,8\,10\,12\,4\,7,\quad   \pi=1\,2\,6\,5\,13\,11,\quad \text{and} \quad \alpha=1\,9\,2\,6\,3\,5\,13\,8\,10\,12\,11\,4\,7,\]  where $m=7$, $r=2$ , $n=6$, $s=2$ and $k=5$. The elements of $\pi$ in $\alpha^{(i)}$ are in boldface to distinguish them from the elements of $\sigma$.  The pairs of partitions $(\lambda,\mu)=\Phi(\alpha)$ can be constructed as follows:

\begin{center}
  \begin{tabular}{c|clcccc}

$i$  & $\alpha^{(i)}$ & $d_i(\pi)$
& ${\rm maj}(\alpha^{(i-1)})-{\rm maj}(\alpha^{(i)})$ & $t(i)$ & ${\rm {\rm des}}(\alpha^{(i-1)})-{\rm des}(\alpha^{(i)})$  \\\hline
6   & 9\,3\,8\,10\,12\, 4\,7 & 0 & 6 & 6 & 1 \\[2pt]
5   &  9\,3\,8\,10\,12\,{\bf 11}\,4\,7  & 1 & 5 & 4 & 1 \\[2pt]
4   & 9\,3\,{\bf 13}\,8\,10\,12\,{\bf 11}\,4\,7& 1 & 3 & 2 & 0 \\[2pt]
3   & 9\,3\,{\bf 5}\,{\bf 13}\,8\,10\,12\,{\bf 11}\,4\,7 & 2 & 5 & 3 & 1 \\[2pt]
2   & 9\,{\bf 6}\,3\,{\bf 5}\,{\bf 13}\,8\,10\,12\,{\bf 11}\,4\,7 & 2 & 4 & 2 & 0  \\[2pt]
1   &  9\,{\bf 2}\,{\bf 6}\,3\,{\bf 5}\,{\bf 13}\,8\,10\,12\,{\bf 11}\,4\,7 & 2 & 5 & 3 & 0   \\[2pt]
0    & {\bf 1}\,9\,{\bf 2}\,{\bf 6}\,3\,{\bf 5}\,{\bf 13}\,8\,10\,12\,{\bf 11}\,4\,7 &   &  &   &     \\[2pt]
\end{tabular}
\end{center}
Hence $\lambda=(6,4,3)$ and $\mu=(3,2,2)$.

In order to prove that the map $\Phi$ is a bijection, we shall reformulate the insertion lemma of Haglund, Loehr and Remmel.  To this end, we first recall the notation of the major increment sequence introduced by Novick \cite{Novick-2010}. Let $\sigma=\sigma_1\cdots \sigma_n\in \mathfrak{S}_n$ and $r\not \in \sigma$. Recall that $\sigma^{(i)}(r)$ denotes the permutation obtained by inserting $r$ before $\sigma_{i+1}$ with the convention that $\sigma_{n+1}=0$. For $0\leq i\leq n$, define the  major increment~
\[{\rm im}(\sigma,i,r)={\rm maj}(\sigma^{(i)}(r))-{\rm maj}(\sigma)\]
and the  major increment sequence
\[MIS(\sigma,r)=({\rm im}(\sigma,0,r),\ldots,{\rm im}(\sigma,n,r)).\]
Combining Lemma \ref{hag-leo-rem} and Lemma \ref{hag-leo-rem-add}, we obtain the following consequence.

\begin{coro}\label{3.1}
Let $\sigma \in \mathfrak{S}_{n}$ with $k$ descents and $r\not \in \sigma$. Then $MIS(\sigma,r)$  is a shuffling of $k+1,k+2,\ldots,n$ and $k,\ldots,1,0$. In particular, ${\rm im}(\sigma,i,r)$ is either $\min\{{\rm im}(\sigma,0,r),\ldots,{\rm im}(\sigma,i-1,r)\}-1$ or $\max\{{\rm im}(\sigma,0,r),\ldots,{\rm im}(\sigma,i-1,r)\}+1$. More precisely, if ${\rm des}(\sigma^{(i)}(r))={\rm des}(\sigma)+1$, then \[{\rm im}(\sigma,i,r)=\max\{{\rm im}(\sigma,0,r),\ldots,{\rm im}(\sigma,i-1,r)\}+1,\]
otherwise,
\[{\rm im}(\sigma,i,r)=\min\{{\rm im}(\sigma,0,r),\ldots,{\rm im}(\sigma,i-1,r)\}-1.\]
\end{coro}

 For example, let $\sigma=5\,1\,6\,2\,4 \in \mathfrak{S}_5$ and $r=3$, we have ${\rm des}(\sigma)=2$ and ${\rm maj}(\sigma)=4$.

\begin{center}
  \begin{tabular}{c|cccc}

$i$ & $\sigma^{(i)}(3)$ & ${\rm maj}(\sigma^{(i)}(3))$
& ${\rm im}(\sigma,i,3)$ & ${\rm des}(\sigma^{(i)}(3))-{\rm des}(\sigma)$\\\hline
0 & {\bf3}\,5\,1\,6\,2\,4 & 6 & 2 &0\\[2pt]
1 & 5\,{\bf 3}\,1\,6\,2\,4 & 7 & 3 &1\\[2pt]
2 & 5\,1\,{\bf 3}\,6\,2\,4 & 5 & 1 &0\\[2pt]
3 & 5\,1\,6\,{\bf 3}\,2\,4 & 8 & 4 &1\\[2pt]
4 & 5\,1\,6\,2\,{\bf 3}\,4 & 4 & 0 &0\\[2pt]
5 & 5\,1\,6\,2\,4\,{\bf 3} & 9 & 5  &1
\end{tabular}
\end{center}
Moreover,  $MIS(\sigma,r)=(2,3,1,4,0,5)$ which is a shuffle of $3,4,5$ and $2,1,0$.

Let $MI{S}_i(\sigma,r)=({\rm im}(\sigma,0,r),\ldots,{\rm im}(\sigma,i-1,r))$ be the first $i$ elements of $MIS(\sigma,r)$. Employing  the insertion lemma of Haglund, Loehr and Remmel,  Novick \cite{Novick-2010} found the following interesting proposition about $MI{S}_i(\sigma,r)$. It turns out that this proposition plays an important role in the proof that  the map $\Phi$ is a bijection. For completeness,  we provide an alternative proof of this proposition with the aid of Corollary  \ref{3.1}. Here we use the common notation $\chi(T)=1$ if the statement $T$ is true and $\chi(T)=0$   otherwise.

\begin{prop}[Novick]\label{MIS}
Suppose that $\sigma$ is  a permutation of length $m$ with $r$ descents. Let $p,q\notin \sigma$ and let $\sigma^{(i-1)}(p)$ be the permutation by inserting $p$ before $\sigma_i$. Then \break $MI{S}_i(\sigma^{(i-1)}(p),q)$ is a permutation of the set $\{{\rm im}(\sigma,j,p)+\chi(q>p)|0\leq j< i\}$.
\end{prop}
\proof
Let $\sigma[i]=\sigma_1\sigma_2\cdots\sigma_i$ be the permutation of first $i$ elements of $\sigma$. By Corollary \ref{3.1}, we find that
$MIS(\sigma[i],p)$ is a permutation of the set $\{0,1,\ldots, i\}$.
Note that ${\rm im}(\sigma[i],i,p)=i\chi(\sigma_i>p)$, so $MI{S}_i(\sigma[i],p)$ is a permutation of the set $\{j-1+\chi(p>\sigma_i)|0<j\le i\}$. By the definition of descents, we find that
\[{\rm im}(\sigma[i],j,p)+d_i(\sigma)={\rm im}(\sigma,j,p) \quad \text{ for} \quad 0\leq j< i.\]
Hence  $MI{S}_i(\sigma,p)$ is a permutation of the set $\{j-1+\chi(p>\sigma_i)+d_i(\sigma)|0<j\le i\}$. Using the same argument, we derive that $MI{S}_i(\sigma^{(i-1)}(p),q)$ is a permutation of the set $\{j-1+\chi(q>p)+d_i(\sigma^{(i-1)}(p))|0<j\le i\}$, where $\sigma^{(i-1)}_i(p)=p$. Note that $d_i(\sigma^{(i-1)}(p))=\chi(p>\sigma_i)+d_i(\sigma)$, so we conclude that  $MI{S}_i(\sigma^{(i-1)}(p),q)$ is a permutation of the set $\{j-1+\chi(q>p)+\chi(p>\sigma_i)+d_i(\sigma)|0<j\le i\}$. Then the proposition follows immediately by comparing  $MI{S}_i(\sigma^{(i-1)}(p),q)$ with $MI{S}_i(\sigma,p)$. This completes the proof. \qed

For example, let $\sigma= 5\,8\,1\,4\,6\,2 \in \mathfrak{S}_6$, $p=7$, $q=9$ and $i=5$. Note that $\sigma^{(4)}(7)=5\,8\,1\,4\,7\,6\,2 \in \mathfrak{S}_7$, and it can be evaluated that   \[MIS(\sigma,7)=(3,\,2,\,4,\,5,\,6,\,1,\,0)\]
and
\[MIS(\sigma^{(4)}(7),9)=(4,\,5,\,3,\,6,\,7,\,2,\,1,\,0).\]
We find that $MI{S}_5(\sigma,7)$ is a permutation of $\{2,\,3,\,4,\,5,\,6\}$  while
$MI{S}_5(\sigma^{(4)}(7),9)$ is a permutation of $\{3,\,4,\,5,\,6,\,7\}$. On the other hand,  $\sigma^{(4)}(9)=5\,8\,1\,4\,9\,6\,2 \in \mathfrak{S}_7$, and it can be checked that
\[MIS(\sigma,9)=(3,\,4,\,2,\,5,\,6,\,1,\,0),\]
and
\[MIS(\sigma^{(4)}(9),7)=(4,\,3,\,5,\,6,\,2,\,7,\,1,\,0).\]
It is clear that   both $MI{S}_5(\sigma,9)$ and $MI{S}_5(\sigma^{(4)}(9),7)$ are  permutations of $\{2,\,3,\,4,\,5,\,6\}$.

With   Corollary \ref{3.1} and Proposition  \ref{MIS} in hand, we are now in a position to show that the map $\Phi$  in Definition \ref{defi-map} is a map from $\mathfrak{S}(\sigma, \pi|k)$ to $ \mathcal{P}_{k-r}(k-s, m)\times \mathcal{P}_{n-k+r}(0,\,k-s)$.

\begin{lem} \label{thm-lem} Suppose that $\sigma\in \mathfrak{S}_m$ and $\pi\in \mathfrak{S}_n$ are two disjoint permutations, where  ${\rm des}(\sigma)=r$ and ${\rm des}(\pi)=s$. Let $\alpha\in \mathfrak{S}(\sigma, \pi|k)$ and $(\lambda,\mu)=\Phi(\alpha)$. Then $\lambda \in \mathcal{P}_{k-r}(k-s,m)$ and $\mu \in \mathcal{P}_{n-k+r}(0,\,k-s)$. Moreover, $ {\rm maj}(\alpha)=|\lambda|+|\mu|+{\rm maj}(\sigma)+{\rm maj}(\pi).$
\end{lem}

\proof Recall that $\sigma=\sigma_1\cdots\sigma_m$ is a permutation   with $r$  descents and  $\pi=\pi_1\cdots\pi_n$ be a permutation  with $s$  descents.  Assume that $\alpha=\alpha_1\cdots \alpha_{n+m}$ is  the shuffle of $\sigma$ and $\pi$ with $k$  descents  with the convention that $\alpha_{n+m+1}=+\infty$. For $1\leq i\leq n$, let $\alpha^{(i)}$ denote the permutation obtained by removing $\pi_1,\pi_{2},\ldots,\pi_{i}$ from $\alpha$ and let $k_i$ be  the position at which $\pi_i$ is inserted into $\alpha^{(i)}$ to yield $\alpha^{(i-1)}$. More precisely, $\alpha^{(i-1)}$ is obtained by inserting $\pi_i$ into $\alpha^{(i)}$ before $\alpha^{(i)}_{k_i}$.  Since $\alpha$ is a shuffle of $\sigma$ and $\pi$, we deduce that $k_1\leq k_2\leq \cdots \leq k_n.$ For $1\leq i\leq n$, define
\begin{equation}\label{defi-{thm-lem}}
T^{(i)}=({\rm im}(\alpha^{(i)},0,\pi_i)-d_i(\pi),\ldots,{\rm im}(\alpha^{(i)},k_i-1,\pi_i)-d_i(\pi)).
\end{equation}
From the definition \eqref{defi-t} of $t(i)$, it's clear that $t(i)$ is the final element of $T^{(i)}$.
By Corollary \ref{3.1}, we see that the elements of $T^{(i)}$ are distinct. So we may assume that $T^{(i)}$ is a permutation of a set $ST^{(i)}$.  In terms of Proposition \ref{MIS}, Novick \cite{Novick-2010} showed that
\begin{equation}\label{prop-novick}
  ST^{(1)} \subseteq ST^{(2)} \subseteq \cdots \subseteq ST^{(n)} \subseteq \{0,1,\ldots, m\}.
\end{equation}
It implies that $0\le t(l) \le m$ for any $1\le l\le n$.

Recall that $\alpha^{(i_1)},\ldots,\alpha^{(i_{k-r})}$ are $k-r$ permutations such that
\[{\rm des}(\alpha^{(i_l)})={\rm des}(\alpha^{(i_l+1)})+1 \quad \text{for} \quad 1\le l\le k-r,\]
where $1\leq i_1<i_2<\cdots <i_{k-r}\leq n$ and  $\alpha^{(j_1)},\ldots,\alpha^{(j_{n-k+r})}$ are  permutations such that
\[{\rm des}(\alpha^{(j_l-1)})={\rm des}(\alpha^{(j_l)}) \quad  \text{for}  \quad 1\le l\le n-k+r,\]
where $1\leq j_1<j_2<\cdots <j_{n-k+r}\leq n$. From Corollary \ref{3.1}, we see that   if ${\rm des}(\sigma^{(i)}(r))={\rm des}(\sigma)+1$, then ${\rm im}(\sigma,i,r)=\max\{{\rm im}(\sigma,0,r),\ldots,{\rm im}(\sigma,i-1,r)\}+1$, otherwise, ${\rm im}(\sigma,i,r)=\min\{{\rm im}(\sigma,0,r),\ldots,{\rm im}(\sigma,i-1,r)\}-1$.   Since $t(i)$ is the final element of $T_i$, it follows that
 $t(i)$   is either the largest element of $T_i$ or the smallest element of $T_i$.  Hence, by \eqref{prop-novick},  we derive that
 \begin{equation}\label{pf-c1}
   m\ge t(i_{k-r}) \ge \cdots \ge t(i_{1}) \ge t(j_1) \ge \cdots\ge t(j_{n-k+r})\geq 0.
 \end{equation}
To prove that $\lambda \in  \mathcal{P}_{k-r}(k-s, m)$ and $\mu \in \mathcal{P}_{n-k+r}(0,\,k-s)$, by \eqref{defi-lambda} and \eqref{defi-mu},  it suffices to show that $t(j_1) \le k-s \le t(i_1)$.
By the definition of $i_1$, we see that  ${\rm des}(\alpha^{(i_1-1)})={\rm des}(\alpha)=k$ and  ${\rm des}(\alpha^{(i_1)})=k-1$, and so  ${\rm maj}(\alpha^{(i_1-1)})-{\rm maj}(\alpha^{(i_1)}) \geq k$. But $d_{i_1}(\pi)\leq {\rm des}(\pi)=s$, it follows that
\begin{equation}\label{pf-c2}
t(i_1)={\rm maj}(\alpha^{(i_1-1)})-{\rm maj}(\alpha^{(i_1)})-d_{i_1}(\pi)\ge k-s.
\end{equation}
Combining \eqref{pf-c1} and \eqref{pf-c2}, we conclude that   $   \lambda=(t(i_{k-r}),t(i_{k-r-1}),\ldots, t(i_{1}))$ is a partition in
$\mathcal{P}_{k-r}(k-s, m)$. From the definition of $j_1$, we see that  ${\rm des}(\alpha^{(j_1-1)})={\rm des}(\alpha^{({j_1})})=k-j_1+1$, and so  ${\rm maj}(\alpha^{(j_1-1)})-{\rm maj}(\alpha^{(j_1)})\le k-j_1+1$. Since $d_{j_1}(\pi)\geq {\rm des}(\pi)-j_1+1=s-j_1+1$, we have
\begin{equation}\label{pf-c3}
t(j_1)={\rm maj}(\alpha^{(j_1-1)})-{\rm maj}(\alpha^{(j_1)})-d_{j_1}(\pi)\le k-s.
\end{equation}
Combining \eqref{pf-c1} and \eqref{pf-c3}, we arrive at  $\mu=(t(j_{1}),t(j_{2}),\ldots, t(j_{n-k+r}))$ is a partition in
$\mathcal{P}_{n-k+r}(0,\,k-s)$.  Moreover,  it is evident from \eqref{defi-t}, \eqref{defi-lambda} and \eqref{defi-mu} that  ${\rm maj}(\alpha)-{\rm maj}(\sigma)={\rm maj}(\alpha^{(0)})-{\rm maj}(\alpha^{(n)})=\sum_{i=1}^n t(i) +{\rm maj}(\pi)=|\lambda|+|\mu|+{\rm maj}(\pi)$.   This completes the proof. \qed

 To prove Theorem \ref{Stanley-c}, we also need to define the inverse map of $\Phi$.

\begin{defi}[The map $\Psi$]\label{defi-psi}  Assume that $\sigma$ and $\pi$ are given  in Theorem  \ref{Stanley-c}.
 Let $\lambda=(\lambda_1,\ldots, \lambda_{k-r}) \in  \mathcal{P}_{k-r}(k-s, m)$ , $\mu=(\mu_1,\ldots,\mu_{n-k+r}) \in \mathcal{P}_{n-k+r}(0,\,k-s)$ and
\begin{equation}\label{part-cond-multiset}
 M^{(n)}=\{\lambda_1,\ldots, \lambda_{k-r}, \mu_1,\ldots, \mu_{n-k+r}\}
\end{equation}
be a multiset consisting  of all parts of  $\lambda$ and $\mu$. All of the elements in $ M^{(n)}$ are listed in   non-increasing order.
The map $\Psi\colon (\lambda,\mu)\rightarrow \alpha$  is defined as follows: Assume that $\alpha^{(n)}=\sigma$ and $k_{n+1}=m+1$. Set $b=0$ and carry out the following procedure.

 {\rm (A)}  Define
 \begin{equation}\label{pf-inverse-a}
  T^{(n-b)}=({\rm im}(\alpha^{(n-b)},0,\pi_{n-b})-d_{n-b}(\pi),\ldots,{\rm im}(\alpha^{(n-b)},k_{n-b+1}-1,\pi_{n-b})-d_{n-b}(\pi) ).
 \end{equation}
 We use $T^{(n-b)}_i$ to denote the $i$-th element of $T^{(n-b)}$. Let $k_{n-b}$ be the largest positive integer such that $T^{(n-b)}_{k_{n-b}}\in M^{(n-b)}$ {\rm(}the existence of $k_{n-b}$ will be proved in Lemma \ref{lem2}{\rm)}.
  Let $\alpha^{(n-b-1)}$ be the permutation obtained by inserting $\pi_{n-b}$ into $\alpha^{(n-b)}$ before  $\alpha_{k_{n-b}}^{(n-b)}$.
    Define
    \begin{equation}\label{pf-inverse-bb}
    M^{(n-b-1)}=M^{(n-b)}\setminus \{T^{(n-b)}_{k_{n-b}}\},
    \end{equation}
  which is a multiset of length $n-b-1$.

{\rm (B)} Replace $b$ by $b+1$. If
$b=n$, then we are done. Otherwise, go back to {\rm (A)}.

Then $\Psi(\lambda,\mu)=\alpha^{(0)}$.
\end{defi}

 For example, let $\sigma=9\,3\,8\,10\,12\,4\,7 \in \mathfrak{S}_7$, $\pi=1\,2\,6\,5\,13\,11\in \mathfrak{S}_6$, where ${\rm des}(\sigma)=2$ and ${\rm des}(\pi)=2$. Given $k=5$, $\lambda=(6,4,3)$ and $\mu=(3,2,2)$, we will recover the shuffle $\alpha$ of $\sigma$ and $\pi$ as follows.  The elements of $\pi$ in $\alpha^{(i)}$ are in boldface to distinguish them from the elements of $\sigma$.

\begin{center}
  \begin{tabular}{c|cllllc}
 $i$ & $\pi_{i}$ & $T^{(i)}$ & $M^{(i)}$ & $k_i$ & $\alpha^{(i)}$\\ \hline  6 & 11 & (3,\,2,\,4,\,5,\,1,\,{\bf 6},\,7,\,0) &\{6,\,4,\,3,\,3,\,2,\,2\}&6& 9\,3\,8\,10\,12\,4\,7  \\[2pt]
 5 & 13 & (3,\,2,\,{\bf 4},\,5,\,6,\,1,\ldots) &\{4,\,3,\,3,\,2,\,2\} &3& 9\,3\,8\,10\,12\,{\bf 11}\,4\,7    \\[2pt]
 4 & 5 & (3,\,4,\,{\bf 2},\,\ldots) &\{3,\,3,\,2,\,2\} &3& 9\,3\,{\bf 13}\,8\,10\,12\, {\bf 11}\,4\,7 \\[2pt]
 3 & 6 & (2,\,{\bf 3},\,4,\,\ldots) &\{3,\,3,\,2\} &2 & 9\,3\,{\bf 5}\,{\bf 13}\,8\,10\,12\, {\bf 11}\,4\,7    \\[2pt]
 2 & 2 & (3,\,{\bf 2},\,\ldots)&\{3,\,2\}&2 & 9\,{\bf 6}\,3\,{\bf 5}\,{\bf 13}\,8\,10\,12\,{\bf 11}\,4\,7     \\[2pt]
 1 & 1 & ({\bf 3},\,2,\,\ldots) &\{3\} &1& 9\,{\bf 2}\,{\bf 6}\,3\,{\bf 5}\,{\bf 13}\,8\,10\,12\,{\bf 11}\,4\,7  \\[2pt]
  0 &   &   &$\emptyset$ & & {\bf 1}\,9\,{\bf 2}\,{\bf 6}\,3\,{\bf 5}\,{\bf 13}\,8\,10\,12\,{\bf 11}\,4\,7
\end{tabular}
\end{center}

Hence $\alpha=\alpha^{(0)}=1\,9\,2\,6\,3\,5\,13\,8\,10\,12\,11\,4\,7$.

We proceed to prove that the map $\Psi$ defined in Definition \ref{defi-psi}  is a map from $\mathcal{P}_{k-r}(k-s, m) \times \mathcal{P}_{n-k+r}(0,\,k-s)$ to $\mathfrak{S}(\sigma, \pi|k)$.

\begin{lem} \label{lem2} Assume that $\sigma$ and $\pi$ are given  in Theorem  \ref{Stanley-c}. Let  $(\lambda,\mu) \in  \mathcal{P}_{k-r}(k-s, m) \times \mathcal{P}_{n-k+r}(0,\,k-s)$ and let  $\alpha=\Psi(\lambda,\mu)$. Then $\alpha \in \mathfrak{S}(\sigma, \pi|k)$.
\end{lem}

\pf  We first show that $\alpha$ is a shuffle of $\sigma$ and $\pi$. To this end, we need to show that $k_{n-b}$   in Definition \ref{defi-psi} exists and $k_{n-b}\leq k_{n-b+1}$ for $0\leq b\leq n-1$.

Using Corollary \ref{3.1}, we find that when $b=0$, $T^{(n)}$ is a permutation of the set $\{0,1,\ldots, m\}$ since $d_{n}(\pi)=0$. From \eqref{part-cond-multiset}, we see that all elements in $M^{(n)}$ are in $T^{(n)}$, and so $k_n$ exists. For $1\leq b\leq n-1$, assume that $k_{n-b+1}$ exists, we proceed to show that $k_{n-b}$ exists.  In light of Proposition \ref{MIS}, we derive that the elements in $T^{(n-b)}$ are the same as the first $k_{n-b+1}$ elements in $T^{(n-b+1)}$. Since $k_{n-b+1}$ is the largest integer such that $T_{k_{n-b+1}}^{(n-b+1)} \in M^{(n-b+1)}$, we deduce that  all of the elements in $M^{(n-b)}$ are located to the left of $T_{k_{n-b+1}}^{(n-b+1)}$ in $T^{(n-b+1)}$, and so   all of the elements in  $M^{(n-b)}$ are also in  $T^{(n-b)}$. It follows that $k_{n-b}$ exists. Moreover,  by definition,  it is easy to see that
$ k_{n-b}\leq k_{n-b+1}.$  Hence we have proven that $\alpha$ is a shuffle of $\sigma$ and $\pi$.

We next show that there are $k$ descents in $\alpha$. Suppose to the contrary that ${\rm des}(\alpha)\neq k$. Assume that ${\rm des}(\alpha)=l<k$.  Let $(\overline{\lambda},\overline{\mu})=\Phi(\alpha)$,   by Lemma \ref{thm-lem}, we derive that $\overline{\lambda} \in \mathcal{P}_{l-r}(l-s,m)$ and $\overline{\mu} \in \mathcal{P}_{ n-l+r}(0, l-s)$, that is,
 \[  m\geq \overline{\lambda}_1\geq \cdots \ge \overline{\lambda}_{l-r}\geq l-s \geq \overline{\mu}_1\geq \cdots \geq \overline{\mu}_{n-l+r}\geq 0. \]
 Let
 \[\overline{M}^{(n)}=\{\overline{\lambda}_1,\ldots,  \overline{\lambda}_{l-r},\overline{\mu}_1, \ldots   \overline{\mu}_{n-l+r}\}\]
 be a multiset consisting of all parts of $\overline{\lambda}$ and $\overline{\mu}$.  All of the elements in $\overline{M}^{(n)}$ are listed in   non-increasing order.
 From the definitions of $\Phi$ and $\Psi$, it is easy to see that
 $\overline{M}^{(n)}$ equals  ${M}^{(n)}$ defined in \eqref{part-cond-multiset}.
 Since $k>l$, we derive that $\overline{\mu}_{k-l}=\lambda_{k-r}\geq k-s$ which contradicts the fact that $\overline{\mu}_{k-l}\leq l-s<k-s$. Applying the same argument, we deduce that   ${\rm des}(\alpha)=l>k$, which is also impossible. Hence we arrive at the conclusion that  ${\rm des}(\alpha)=k$.
Therefore, $\Psi$ is a map from  $ \mathcal{P}_{k-r}(k-s, m)\times \mathcal{P}_{n-k+r}(0,\,k-s)$ to $\mathfrak{S}(\sigma,\pi|k)$. This completes the proof. \qed

The final part of this section is to provide a proof of Theorem \ref{Stanley-c} using Lemma \ref{thm-lem} and Lemma \ref{lem2}.

\noindent{\it Proof of Theorem \ref{Stanley-c}:}  Let $\alpha \in \mathfrak{S}(\sigma, \pi|k)$. Utilizing Lemma \ref{thm-lem}, we find that $\Phi(\alpha)$ belongs to $\mathcal{P}_{k-r}(k-s, m)\times \mathcal{P}_{n-k+r}(0,\,k-s)$. Combining the definition of $\Phi$ and the definition of $\Psi$, we deduce that $\Psi(\Phi(\alpha))=\alpha$.

Conversely, let  $(\lambda,\mu)\in \mathcal{P}_{k-r}(k-s, m)\times \mathcal{P}_{n-k+r}(0,\,k-s)$. Invoking Lemma \ref{lem2}, we know that $\Psi(\lambda,\mu)\in \mathfrak{S}(\sigma, \pi|k)$. By virtue of Definition \ref{defi-map} and Definition \ref{defi-psi}, we obtain that  $\Phi(\Psi(\lambda,\mu))=(\lambda,\mu)$.

Therefore, the map $\Phi$ is a bijection between $\mathfrak{S}(\sigma, \pi|k)$ and  $\mathcal{P}_{k-r}(k-s, m)\times \mathcal{P}_{n-k+r}(0,\,k-s)$. This completes the proof. \qed

\section{Refinements}

In this section, we first state four refinements of Stanley's shuffle theorem and then provide proofs of these refinements with the help of  Lemma \ref{thm-lem} and Lemma \ref{lem2}.

Suppose that $\sigma=\sigma_1\cdots\sigma_{m}\in \mathfrak{S}_m$
and $\pi=\pi_1\cdots \pi_n\in \mathfrak{S}_n$ are two disjoint permutations, where ${\rm des}(\sigma)=r$ and ${\rm des}(\pi)=s$. From Theorem 3.1, we see that  Stanley's shuffle theorem  is equivalent to the statement that there is a bijection $\Phi$ between $\mathfrak{S}(\sigma,\pi|k)$ and $\mathcal{P}_{k-r}(k-s, m) \times \mathcal{P}_{n-k+r}(0,\,k-s)$ such that for    $\alpha \in \mathfrak{S}(\sigma,\pi|k)$, we have $(\lambda,\mu)=\Phi(\alpha)\in  \mathcal{P}_{k-r}(k-s, m) \times \mathcal{P}_{n-k+r}(0,\,k-s)$, namely,
\[  m\geq \lambda_1\geq \cdots \ge \lambda_{k-r}\geq k-s \geq  {\mu}_1\geq \cdots \geq {\mu}_{n-k+r}\geq 0. \]

 In our first refinement,  we restrict our attention to   the
   subset of $\mathfrak{S}(\sigma,\pi|k)$,
   denoted $\mathfrak{S}^{sb}(\sigma,\pi|k)$
   such that for  $\alpha \in \mathfrak{S}^{sb}
   (\sigma,\pi|k)$,  we have
   $(\lambda,\mu)=\Phi(\alpha)$ where $\mu_{n-
   k+r}\geq 1$. This is also the reason that
   we denote this subset by $\mathfrak{S}^{sb}
   (\sigma,\pi|k)$. The images of $\alpha \in
   \mathfrak{S}^{sb}(\sigma,\pi|k)$ under the
   bijection $\Phi$ give more restrictions on  the
   smallest part of the partition $\mu$.

In the same vein, our second refinement is defined on the subset  $\mathfrak{S}^{la}(\sigma,\pi|k)$  of $\mathfrak{S}(\sigma,\pi|k)$ and the  third refinement is defined on the subset  $\mathfrak{S}^{sa}(\sigma,\pi|k)$. For  $\alpha \in \mathfrak{S}^{la}(\sigma,\pi|k)$,  we have $(\lambda,\mu)=\Phi(\alpha)$ satisfying $\lambda_1=m$, whereas,   for  $\alpha \in \mathfrak{S}^{sa}(\sigma,\pi|k)$,  we have $(\lambda,\mu)=\Phi(\alpha)$ such that $\lambda_{k-r}=k-s$. The fourth refinement involves  the subset  $\mathfrak{S}^{lb}(\sigma,\pi|k)$ such that for  $\alpha \in \mathfrak{S}^{lb}(\sigma,\pi|k)$,  we have $(\lambda,\mu)=\Phi(\alpha)$ satisfying $\mu_{1}=k-s$.

To express these four refinements in terms of generating function, we need to consider the generating functions of the following two special sets of partitions.
  Let $\mathcal{P}_n(t,=m)$ denote  the set of partitions $\lambda$ such that $\ell(\lambda)=n$, $\lambda_n\geq t$ and $\lambda_1=m$. Using \eqref{int-GassCoef} as a starting point, it's not difficult to  derive that
\begin{equation}\label{int-GassCoeft-2}
q^{(n-1)t+m}{n+m-t-1 \brack n-1}=\sum_{\lambda \in \mathcal{P}_n(t,=m)} q^{|\lambda|}.
\end{equation}

Similarly, let $\mathcal{P}_n(=t,m)$ denote
 the set of partitions $\lambda$ such that $\ell(\lambda)=n$, $\lambda_n=t$ and $\lambda_1\leq m$, we have
\begin{equation}\label{int-GassCoeft-3}
q^{nt}{n+m-t-1 \brack n-1}=\sum_{\lambda \in \mathcal{P}_n(=t,m)} q^{|\lambda|},
\end{equation}
which is   required in the proof of the third refinement of Stanley's shuffle theorem.

\subsection{Statements of refinements}

Suppose that $\sigma=\sigma_1\cdots\sigma_{m}\in \mathfrak{S}_m$
and $\pi=\pi_1\cdots \pi_n\in \mathfrak{S}_n$ are two disjoint permutations.  Let
$\mathfrak{S}^{sb}(\sigma,\pi)$ denote the set of shuffles
$\alpha=\alpha_1\cdots\alpha_{n+m}$ of $\sigma$ and
$\pi$ such that $\alpha_{n+m}=\min \{\pi_n,\sigma_m\}$.
For example,
\[\mathfrak{S}^{sb}(2\,6\,3, {\bf 1}\,{\bf 4})=\{ 2\,6\,{\bf 1}\,{\bf 4}\,3, 2\, {\bf 1}\,6\,{\bf 4}\,3,  {\bf 1}\, 2\,6\,{\bf 4}\,3, 2\,{\bf 1}\,{\bf 4}\,6\,3, {\bf 1}\,2\,{\bf 4}\,6\,3, {\bf 1}\,{\bf 4}\,2\,6\,3\}
  \]
  and
  \[\mathfrak{S}^{sb}(2\,8\,5, {\bf 1}\,{\bf 4})=\{ 2\,8\,5\,{\bf 1}\,{\bf 4}, 2\,8\,{\bf 1}\,5\,{\bf 4} , 2\,{\bf 1}\,8\,5\,{\bf 4}, {\bf 1}\,2\,8\,5\,{\bf 4}\}.
  \]

Our first refinement of Stanley's shuffle theorem is the following.

\begin{thm}[The first refinement] \label{stanley-r2}
Let $\sigma=\sigma_1\cdots\sigma_{m}\in \mathfrak{S}_m$ and $\pi=\pi_1\cdots \pi_n\in \mathfrak{S}_n$ be two disjoint permutations.
\begin{itemize}
\item If $\pi_n>\sigma_m$, then
 \begin{align}
    \sum_{\alpha\in \mathfrak{S}^{sb}(\sigma,\pi) \atop {\rm des}(\alpha)=k}q^{{\rm maj}(\alpha)}&=
    {\ell(\sigma)-{\rm des}(\sigma)+{\rm des}(\pi) \brack k-{\rm des}(\sigma)}{\ell(\pi)-{\rm des}(\pi)+{\rm des}(\sigma)-1 \brack k-{\rm des}(\pi)-1} \nonumber\\[5pt]
    &\quad \quad \quad \times q^{{\rm maj}(\sigma)+{\rm maj}(\pi)+\ell(\pi)+(k-{\rm des}(\pi)-1)(k-{\rm des}(\sigma))}.
\end{align}
\item If $\pi_n<\sigma_m$, then
\begin{align}
    \sum_{\alpha\in \mathfrak{S}^{sb}(\pi,\sigma) \atop {\rm des}(\alpha)=k}q^{{\rm maj}(\alpha)}&={\ell(\sigma)-{\rm des}(\sigma)+{\rm des}(\pi)-1 \brack k-{\rm des}(\sigma)-1}{\ell(\pi)-{\rm des}(\pi)+{\rm des}(\sigma) \brack k-{\rm des}(\pi)}\nonumber\\[5pt]
    &\quad \quad \quad \times q^{{\rm maj}(\sigma)+{\rm maj}(\pi)+\ell(\sigma)+(k-{\rm des}(\pi))(k-{\rm des}(\sigma)-1)}. \label{sym-m}
\end{align}
\end{itemize}
\end{thm}
As an immediate corollary of Theorem \ref{stanley-r2}, we obtain the following result.
\begin{coro}
Suppose that $\sigma=\sigma_1\cdots\sigma_{m}\in \mathfrak{S}_m$ and $\pi=\pi_1\cdots \pi_n\in \mathfrak{S}_n$ are two disjoint permutations. Let $\mathfrak{S}^{sb}(\sigma,\pi|k)$ denote the set of shuffles $\alpha=\alpha_1\cdots \alpha_{n+m}$ of $\sigma$ and $\pi$ such that $\alpha_{n+m}=\min\{\sigma_{m}, \pi_n\}$ with $k$ descents.
\begin{itemize}
\item If $\pi_n>\sigma_m$, then
\[\# \mathfrak{S}^{sb}(\sigma,\pi|k)={\ell(\sigma)-{\rm des}(\sigma)+{\rm des}(\pi) \choose k-{\rm des}(\sigma)}{\ell(\pi)-{\rm des}(\pi)+{\rm des}(\sigma)-1 \choose k-{\rm des}(\pi)-1}.\]
\item If $\pi_n<\sigma_m$, then
\[\# \mathfrak{S}^{sb}(\sigma,\pi|k)={\ell(\sigma)-{\rm des}(\sigma)+{\rm des}(\pi)-1 \choose k-{\rm des}(\sigma)-1}{\ell(\pi)-{\rm des}(\pi)+{\rm des}(\sigma) \choose k-{\rm des}(\pi)}.\]
\end{itemize}
\end{coro}

We now turn to state the second refinement of Stanley's shuffle theorem. In doing so, we need to introduce a new shuffle subset.

\begin{defi}\label{defi-refi-1} Let $\sigma=\sigma_1\cdots\sigma_{m}\in \mathfrak{S}_m$ and $\pi=\pi_1\cdots \pi_n\in \mathfrak{S}_n$ be two disjoint permutations and let $\alpha$ be a shuffle of $\sigma$ and $\pi$.  If  $\pi_n>\sigma_m$, then set $\gamma=\pi$ and $\delta=\sigma$. If $\pi_n<\sigma_m$, then set $\gamma=\sigma$ and $\delta=\pi$. Let $a_i$ be the largest number such that  $\delta_{a_i}$ is before $\gamma_i$ in $\alpha$.  Then $\alpha$ can be represented in the following form:
\begin{equation}\label{defi-sec}
\alpha=\delta_1\cdots \delta_{a_1}\gamma_1\delta_{a_1+1}\cdots \delta_{a_2}\gamma_2\delta_{a_2+1}\cdots \delta_{a_{\ell(\gamma)}} \gamma_{\ell(\gamma)} \delta_{a_{\ell(\gamma)}+1}\cdots \delta_{\ell(\delta)}.
\end{equation}
Here we assume that $\delta_{0}=0$, $\delta_{a_{\ell(\gamma)+1}}=0$ and $\gamma_{\ell(\gamma)+1}=+\infty$.  Let ${\rm ld}(\gamma)$ denote the largest descent of $\gamma$   if it exists and $0$ otherwise.
A shuffle $\alpha$ of $\sigma$ and $\pi$ is defined to be in the  set   $\mathfrak{S}^{la}(\sigma,\pi)$ if there exists ${\rm ld}(\gamma)\leq j\leq \ell(\gamma)$ in \eqref{defi-sec} satisfying the following four conditions{\rm :}
\begin{itemize}
\item[{\rm (1)}] $\gamma_{j+1}>\delta_{a_{j+1}}>\delta_{a_{j+1}+1}>\cdots>\delta_{\ell(\delta)}.$
\item[{\rm (2)}] If $\gamma_j<\delta_{a_{j+1}}$, then $\delta_{a_j}=\delta_{a_{j+1}}$.
\item[{\rm (3)}] If $j={\rm ld}(\gamma)$, then $\delta_{a_j}=\delta_{a_{j+1}}$.
\item[{\rm (4)}] If $\gamma_j>\delta_{a_{j+1}}$ and $j>{\rm ld}(\gamma)$, then $\delta_{a_j}\neq \delta_{a_{j+1}}$ and $\gamma_j>\delta_{a_{j}+1}>\delta_{a_{j}+2}>\cdots>\delta_{a_{j+1}}$. Moreover, either $\delta_{a_j}>\gamma_j$ or $\delta_{a_j}<\delta_{a_j+1}$.
\end{itemize}
\end{defi}

For example, let $\sigma=8\,13\,5\,3\,2\,1 \in \mathfrak{S}_6$ and $\pi=  9 \, 4 \,  7 \,  11 \,  12\in \mathfrak{S}_5$. Since $\sigma_6=1<\pi_5=12$,  set $\gamma=\pi={\bf 9}\,{\bf 4}\,{\bf 7}\,{\bf 11}\,{\bf 12}$ and $\delta=\sigma=8\,13\,5\,3\,2\,1$. Let
\[\begin{array}{cccccccccccccc}
\alpha&=&8&{\bf9}&13&5&\underline{{\bf 4}}&{\bf 7}&3&2&{\bf 11}&{\bf 12}&1.\\
&&&\downarrow&&&\downarrow&\downarrow&&&\downarrow&\downarrow&\\
&&&\gamma_1&&&\gamma_2&\gamma_3&&&\gamma_4&\gamma_5&\\
\end{array}
\]
We see that $\delta_{a_1}=8$, $\delta_{a_2}=\delta_{a_3}=5$,  $\delta_{a_4}=\delta_{a_5}=2$. It can be checked that $\alpha$ fits the conditions (1) and (2) in Definition \ref{defi-refi-1}, where $j=2$. Hence
$\alpha \in \mathfrak{S}^{la}(\sigma,\pi)$.

As other examples, set  $\alpha=8\,  13\,  5\, 3\, \underline{{\bf 9}}\, {\bf 4}\, 2\, {\bf 7}\,{\bf 11}\,1\, {\bf 12}$. It can be checked that it meets the conditions (1) and (3) in Definition \ref{defi-refi-1}, where $j=1$. If we set  $\alpha=8\,  13\, 5\, {\bf 9}\,  \underline{{\bf 4}}\,3\, {\bf 7}\, 2\, {\bf 11}\,1\, {\bf 12}$, then $\alpha$ satisfies the conditions (1) and (4) in Definition \ref{defi-refi-1}, where $j=2$.

 The second refinement of Stanley's shuffle theorem can be stated as follows.

\begin{thm}[The second refinement] \label{stanley-r1}
Let $\sigma=\sigma_1\cdots\sigma_{m}\in \mathfrak{S}_m$
and $\pi=\pi_1\cdots \pi_n\in \mathfrak{S}_n$ be two
disjoint permutations.

\begin{itemize}
\item If $\pi_n>\sigma_m$, then
\begin{align}
    \sum_{\alpha\in \mathfrak{S}^{la}(\sigma,\pi) \atop {\rm des}(\alpha)=k}q^{{\rm maj}(\alpha)}&={\ell(\sigma)-{\rm des}(\sigma)+{\rm des}(\pi)-1 \brack k-{\rm des}(\sigma)-1}{\ell(\pi)-{\rm des}(\pi)+{\rm des}(\sigma) \brack k-{\rm des}(\pi)} \nonumber \\[5pt]
    &\quad \quad \quad \times q^{{\rm maj}(\sigma)+{\rm maj}(\pi)+\ell(\sigma)+(k-{\rm des}(\pi))(k-{\rm des}(\sigma)-1)}.
   \end{align}
\item If $\pi_n<\sigma_m$, then
\begin{align}
    \sum_{\alpha\in \mathfrak{S}^{la}(\sigma,\pi) \atop {\rm des}(\alpha)=k}q^{{\rm maj}(\alpha)}&={\ell(\sigma)-{\rm des}(\sigma)+{\rm des}(\pi)\brack k-{\rm des}(\sigma)}{\ell(\pi)-{\rm des}(\pi)+{\rm des}(\sigma)-1  \brack k-{\rm des}(\pi)-1} \nonumber \\[5pt]
    &\quad \quad \quad \times q^{{\rm maj}(\sigma)+{\rm maj}(\pi)+\ell(\pi)+(k-{\rm des}(\pi)-1)(k-{\rm des}(\sigma))}.
    \end{align}
    \end{itemize}
\end{thm}

The following result is immediate from Theorem \ref{stanley-r1}.

\begin{coro}
Suppose that  $\sigma=\sigma_1\cdots\sigma_{m}\in \mathfrak{S}_m$
and $\pi=\pi_1\cdots \pi_n\in \mathfrak{S}_n$   are two disjoint permutations. Let $\mathfrak{S}^{la}(\sigma,\pi|k)$ denote the set of shuffles of $\sigma$ and $\pi$ with $k$ descents, satisfying the conditions in Definition \ref{defi-refi-1}.
\begin{itemize}
\item If $\pi_n>\sigma_m$, then
\[\# \mathfrak{S}^{la}(\sigma,\pi|k)={\ell(\sigma)-{\rm des}(\sigma)+{\rm des}(\pi)-1 \choose k-{\rm des}(\sigma)-1}{\ell(\pi)-{\rm des}(\pi)+{\rm des}(\sigma) \choose k-{\rm des}(\pi)}.\]
\item If $\pi_n<\sigma_m$, then
\[\# \mathfrak{S}^{la}(\sigma,\pi|k)={\ell(\sigma)-{\rm des}(\sigma)+{\rm des}(\pi) \choose k-{\rm des}(\sigma)}{\ell(\pi)-{\rm des}(\pi)+{\rm des}(\sigma)-1 \choose k-{\rm des}(\pi)-1}.\]
\end{itemize}
\end{coro}

The combination of  Theorem \ref{stanley-r2} and Theorem \ref{stanley-r1} leads to the following result, which seems to be of interest.

\begin{thm} \label{stanley-r1&r2} Suppose that $\pi=\pi_1\cdots \pi_n\in \mathfrak{S}_n$. Let $\sigma=\sigma_1\cdots \sigma_m$
and $\sigma'=\sigma'_1\cdots \sigma'_m$ be two permutations in $\mathfrak{S}_m$ such that ${\rm des}(\sigma)={\rm des}(\sigma')$ and ${\rm maj}(\sigma)={\rm maj}(\sigma')$. If $\sigma_m<\pi_n<\sigma'_m$, then
\begin{equation}
 \sum_{\alpha\in \mathfrak{S}^{sb}(\sigma',\pi) \atop {\rm des}(\alpha)=k}q^{{\rm maj}(\alpha)}=\sum_{\alpha\in \mathfrak{S}^{la}(\sigma,\pi) \atop {\rm des}(\alpha)=k}q^{{\rm maj}(\alpha)}
\end{equation}
and
\begin{equation}
 \sum_{\alpha\in \mathfrak{S}^{sb}(\sigma,\pi) \atop {\rm des}(\alpha)=k}q^{{\rm maj}(\alpha)}=\sum_{\alpha\in \mathfrak{S}^{la}(\sigma',\pi) \atop {\rm des}(\alpha)=k}q^{{\rm maj}(\alpha)}.
\end{equation}
\end{thm}
It would be interesting to give a direct combinatorial proof of Theorem \ref{stanley-r1&r2}.

In order to present the third refinement of Stanley's shuffle theorem, we need to define the following shuffle subset.

\begin{defi}\label{defi-refi-3} Let $\sigma=\sigma_1\cdots\sigma_{m}\in \mathfrak{S}_m$ and $\pi=\pi_1\cdots \pi_n\in \mathfrak{S}_n$ be two disjoint permutations and let $\alpha=\alpha_1\cdots \alpha_{n+m}$ be a shuffle of $\sigma$ and $\pi$ with the convention that $\alpha_{n+m+1}=0$. Suppose that $\pi_j$ is the $p_j$-th element of $\alpha$ (that is, $\alpha_{p_j}=\pi_j$), and there are $q_j$ elements  from $\sigma$ before $\alpha_{p_j}$ in $\alpha$, say $\sigma_1,\sigma_2,\ldots, \sigma_{q_j}$  {\rm (}$q_j$ could be zero and assume that $\sigma_0=0${\rm )}.  Let ${\rm sd}(\pi)$ denote the smallest descent of $\pi$ if it exists and $n$ otherwise. If there exists  $1\leq j\leq {\rm sd}(\pi)$ such that $\sigma_1>\cdots>\sigma_{q_j}>\alpha_{p_j}$ and  either    $\sigma_{q_j}<\alpha_{p_j+1}$ or   $\alpha_{p_j}>\alpha_{p_j+1}$, then we say that the shuffle $\alpha$ is in the  set   $\mathfrak{S}^{sa}(\sigma,\pi)$.
\end{defi}
For example,
\[\mathfrak{S}^{sa}(2\,6\,3,{\bf 1}\,{\bf 4})=\{2\,\underline{ {\bf 1}}\,{\bf 4}\,6\,3, 2\,\underline{ {\bf 1}}\,6\,{\bf 4}\,3, 2\,\underline{ {\bf 1}}\,6\,3\,{\bf 4}, {\bf 1}\,\underline{ {\bf 4}}\,2\,6\,3\},\]
where the element satisfying the conditions in Definition \ref{defi-refi-3} is underlined.

We obtain the following refinement of Stanley's shuffle theorem involving  the set $\mathfrak{S}^{sa}(\sigma,\pi)$.

\begin{thm}[The third refinement] \label{stanley-r3}
Let $\sigma$
and $\pi$ be two
disjoint permutations. Then
 \begin{align}
    \sum_{\alpha\in \mathfrak{S}^{sa}(\sigma,\pi) \atop {\rm des}(\alpha)=k}q^{{\rm maj}(\alpha)}&={\ell(\sigma)-{\rm des}(\sigma)+{\rm des}(\pi)-1 \brack k-{\rm des}(\sigma)-1}{\ell(\pi)-{\rm des}(\pi)+{\rm des}(\sigma) \brack k-{\rm des}(\pi)} \nonumber \\[5pt]
    &\quad \quad \quad \times q^{{\rm maj}(\sigma)+{\rm maj}(\pi)+(k-{\rm des}(\pi))(k-{\rm des}(\sigma))}.
   \end{align}
\end{thm}

Let $q\rightarrow 1$ in Theorem \ref{stanley-r3}, we obtain the following consequence.

\begin{coro}
Suppose that $\sigma$ and $\pi$ are two disjoint permutations. Let $\mathfrak{S}^{sa}(\sigma,\pi|k)$ denote the set of shuffles of $\sigma$ and $\pi$ with $k$ descents satisfying the conditions in Definition \ref{defi-refi-3}. Then
\[\# \mathfrak{S}^{sa}(\sigma,\pi|k)={\ell(\sigma)-{\rm des}(\sigma)+{\rm des}(\pi)-1 \choose k-{\rm des}(\sigma)-1}{\ell(\pi)-{\rm des}(\pi)+{\rm des}(\sigma) \choose k-{\rm des}(\pi)}.\]
\end{coro}

Our last refinement of Stanley's shuffle theorem is defined on the following shuffle subset.

\begin{defi}\label{defi-refi-4} Let $\sigma=\sigma_1\cdots\sigma_{m}\in \mathfrak{S}_m$ and $\pi=\pi_1\cdots \pi_n\in \mathfrak{S}_n$ be two disjoint permutations and let $\alpha=\alpha_1\cdots \alpha_{n+m}$ be a shuffle of $\sigma$ and $\pi$ with the convention that $\alpha_{n+m+1}=0$. Suppose that $\pi_j$ is the $p_j$-th element of $\alpha$ {\rm (}that is, $\alpha_{p_j}=\pi_j${\rm )} and there are $q_j$ elements  from $\sigma$ before $\alpha_{p_j}$ in $\alpha$ {\rm (}$q_j$ could be zero and assume that $\sigma_0=0${\rm )}, say $\sigma_1,\sigma_2,\ldots, \sigma_{q_j}$.  Let ${\rm sa} (\pi)$ denote the smallest ascent of $\pi$  if it exists and $n$ otherwise.
If there exists $1\leq j\leq {\rm sa}(\pi)$ such that $\sigma_1<\cdots<\sigma_{q_j}<\alpha_{p_j}$ and either    $\sigma_{q_j}>\alpha_{p_j+1}$ or  $\alpha_{p_j}<\alpha_{p_j+1}$, then we say that the shuffle $\alpha$ is in the  set   $\mathfrak{S}^{lb}(\sigma,\pi)$.
\end{defi}
For example,
\[\mathfrak{S}^{lb}(2\,6\,3,{\bf 1}\,{\bf 4})=\{\underline{ {\bf 1}}\,{\bf 4}\,2\,6\,3, \underline{ {\bf 1}}\,2\,{\bf 4}\,6\,3, \underline{ {\bf 1}}\,2\,6\,{\bf 4}\,3, \underline{ {\bf 1}}\,2\,6\,3\,{\bf 4}\},\]
Here we underline the element that meets the conditions in Definition \ref{defi-refi-4}.

By considering the set $ \mathfrak{S}^{lb}(\sigma,\pi)$, we can refine Stanley's shuffle theorem as follows.

\begin{thm}[The fourth refinement] \label{stanley-r4}
Let $\sigma$ and $\pi$ be two disjoint permutations. Then
 \begin{align}
    \sum_{\alpha\in \mathfrak{S}^{lb}(\sigma,\pi) \atop {\rm des}(\alpha)=k}q^{{\rm maj}(\alpha)}&={\ell(\sigma)-{\rm des}(\sigma)+{\rm des}(\pi) \brack k-{\rm des}(\sigma)}{\ell(\pi)-{\rm des}(\pi)+{\rm des}(\sigma)-1 \brack k-{\rm des}(\pi)} \nonumber \\[5pt]
    &\quad \quad \quad \times q^{{\rm maj}(\sigma)+{\rm maj}(\pi)+(k-{\rm des}(\pi))(k-{\rm des}(\sigma)+1)}.
   \end{align}
\end{thm}

Using Theorem \ref{stanley-r4}, we obtain the following consequence.

\begin{coro}
Suppose that $\sigma$ and $\pi$ are two disjoint permutations. Let $\mathfrak{S}^{lb}(\sigma,\pi|k)$ denote the set of shuffles of $\sigma$ and $\pi$ with $k$ descents satisfying the conditions in Definition \ref{defi-refi-4}. Then
\[\# \mathfrak{S}^{lb}(\sigma,\pi|k)={\ell(\sigma)-{\rm des}(\sigma)+{\rm des}(\pi) \choose k-{\rm des}(\sigma)}{\ell(\pi)-{\rm des}(\pi)+{\rm des}(\sigma)-1 \choose k-{\rm des}(\pi)}.\]
\end{coro}

Combining Theorem \ref{stanley-r3} and Theorem \ref{stanley-r4}, together with Theorem \ref{Stanley-c},
we arrive at the following consequence. It would be interesting to give a direct combinatorial proof of Corollary \ref{stanley-r3&r4}.

\begin{coro}\label{stanley-r3&r4} Let $\sigma$ and $\pi$ be two disjoint permutations. Then
 \begin{equation*}
     \sum_{\alpha\in \mathfrak{S}^{lb}(\sigma,\pi) \atop {\rm des}(\alpha)=k}q^{{\rm maj}(\alpha)}+\sum_{\alpha\in \mathfrak{S}^{sa}(\pi,\sigma) \atop {\rm des}(\alpha)=k}q^{{\rm maj}(\alpha)}=\sum_{\alpha\in \mathfrak{S}(\sigma,\pi) \atop {\rm des}(\alpha)=k}q^{{\rm maj}(\alpha)}.
 \end{equation*}
\end{coro}

 \subsection{Proofs}

 In this subsection, we prove Theorem \ref{stanley-r2}, Theorem \ref{stanley-r1} , Theorem \ref{stanley-r3} and Theorem \ref{stanley-r4}   by refining the map $\Phi$ given in Definition \ref{defi-map}  on their corresponding shuffle subsets.

Based on \eqref{int-GassCoeft} and using Theorem \ref{Stanley-c}, we see that the proof of Theorem \ref{stanley-r2} is equivalent to the proof of the following combinatorial statement.

 \begin{thm}\label{Stanley-c-r2} Suppose that $\sigma=\sigma_1\cdots\sigma_m \in \mathfrak{S}_m$ and $\pi=\pi_1\cdots \pi_n \in \mathfrak{S}_n$ are two disjoint permutations, where ${\rm des}(\sigma)=r$  and  ${\rm des}(\pi)=s$.
 If $\pi_n>\sigma_m$, then the  map $\Phi$ given in Definition \ref{defi-map} is   a bijection between $\mathfrak{S}^{sb}(\sigma,\pi|k)$ and $ \mathcal{P}_{k-r}(k-s, m)\times \mathcal{P}_{n-k+r}(1,k-s)$. If $\pi_n<\sigma_m$, then the  map $\Phi$ given in Definition \ref{defi-map} is   a bijection between $\mathfrak{S}^{sb}(\sigma,\pi|k)$ and $ \mathcal{P}_{k-s}(k-r, n)\times \mathcal{P}_{m-k+s}(1,k-r)$.
 \end{thm}

\proof Without loss of generality, we may prove this theorem in the case $\pi_n>\sigma_m$. The case  $\pi_n<\sigma_m$ can be justified with the same argument by exchanging $\sigma$ and $\pi$.

Let $\alpha \in \mathfrak{S}^{sb}(\sigma,\pi|k)$. Assume that $\Phi(\alpha)=(\lambda,\mu)$. From Lemma \ref{thm-lem}, we see that $\lambda \in \mathcal{P}_{k-r}(k-s,m)$ and $\mu \in \mathcal{P}_{n-k+r}(0,\,k-s)$, namely,
\[  m\geq \lambda_1\geq \cdots \ge \lambda_{k-r}\geq k-s \geq  {\mu}_1\geq \cdots \geq {\mu}_{n-k+r}\geq 0. \]
We proceed to show that ${\mu}_{n-k+r}\geq 1$ if $\mu\ne \emptyset$. Recall that $T^{(i)}$ is defined as \eqref{defi-{thm-lem}} for $1\leq i\leq n$. Since $\alpha \in \mathfrak{S}^{sb}(\sigma,\pi|k)$, it follows from Corollary \ref{3.1} that $0\notin T^{(n)}$. Using the relation \eqref{prop-novick}, we find that $0\notin T^{(i)}$ for $1\leq i\leq n$. By the definition of $\Phi$, we see that  there exists $i$ such that ${\mu}_{n-k+r}$ is the final element of  $T^{(i)}$, it implies that ${\mu}_{n-k+r}\geq 1$, and so
$\Phi(\alpha)=(\lambda,\mu)\in \mathcal{P}_{k-r}(k-s, m)\times \mathcal{P}_{ n-k+r}(1,k-s)$.

Conversely, let $\lambda \in \mathcal{P}_{k-r}(k-s,m)$ and $\mu \in \mathcal{P}_{ n-k+r}(1, k-s)$. Assume that $\Psi(\lambda,\mu)=\overline{\alpha}$, where the map $\Psi$ is given in Definition \ref{defi-psi}.
In light of  Lemma \ref{lem2}, we derive that $\overline{\alpha}=\overline{\alpha}_1\cdots \overline{\alpha}_{n+m}$  is a shuffle of $\sigma$ and $\pi$ with $k$ descents.
 To prove that  $\overline{\alpha} \in \mathfrak{S}^{sb}(\sigma,\pi)$, it suffices to show that $\overline{\alpha}_{n+m}=\sigma_m$. Suppose to the contrary that $\overline{\alpha}_{n+m}\neq \sigma_m$, that is, $\overline{\alpha}_{n+m}= \pi_n$. Since $\Psi$ and $\Phi$ are inverses of each other, we have $\Phi(\overline{\alpha})=\Phi(\Psi(\lambda,\mu))=(\lambda,\mu)$.  Let $\overline{\alpha}^{(i)}$ denote the permutation obtained by removing $\pi_1,\ldots, \pi_i$ from $\overline{\alpha}$. Clearly, $\overline{\alpha}^{(n)}=\sigma$.   Since $\overline{\alpha}_{n+m}= \pi_n$ and under the condition that $\sigma_m<\pi_n$, we find that ${\rm des}(\overline{\alpha}^{(n-1)})={\rm des}(\overline{\alpha}^{(n)})$,
 and by  the definition of $\Phi$, we deduce that $\mu_{n-k+r}={\rm maj}(\overline{\alpha}^{(n-1)})-{\rm maj}(\overline{\alpha}^{(n)})-d_n(\pi)=0$, which contradicts the condition that $\mu \in \mathcal{P}_{ n-k+r}(1, k-s)$, that is, $\mu_{n-k+r}\geq 1$. Hence  $\overline{\alpha}_{n+m}=\sigma_m$, and so $\overline{\alpha} \in \mathfrak{S}^{sb}(\sigma,\pi)$. This completes the proof.
\qed

We next give a proof of Theorem \ref{stanley-r1}. According to  \eqref{int-GassCoeft} and \eqref{int-GassCoeft-2}, and by Theorem \ref{Stanley-c}, it suffices to prove the following combinatorial statement.

 \begin{thm}\label{Stanley-c-r1} Suppose that $\sigma=\sigma_1\cdots\sigma_m \in \mathfrak{S}_m$ and $\pi=\pi_1\cdots \pi_n \in \mathfrak{S}_n$ are two disjoint permutations, where ${\rm des}(\sigma)=r$  and  ${\rm des}(\pi)=s$. If $\pi_n>\sigma_m$, then the  map $\Phi$ given in Definition \ref{defi-map} is   a bijection between $\mathfrak{S}^{la}(\sigma,\pi|k)$ and $ \mathcal{P}_{k-r}(k-s, =m)\times \mathcal{P}_{n-k+r}(0,\,k-s)$. If $\pi_n<\sigma_m$, then the  map $\Phi$ given in Definition \ref{defi-map} is   a bijection between $\mathfrak{S}^{la}(\sigma,\pi|k)$ and $ \mathcal{P}_{k-s}(k-r, =n)\times \mathcal{P}_{n-k+s}(0,k-r)$.
 \end{thm}

\proof Similar to Theorem \ref{Stanley-c-r2}, it suffices to prove that this theorem  is true when $\pi_n>\sigma_m$. Let $\alpha \in \mathfrak{S}^{la}(\sigma,\pi|k)$. Assume that $\Phi(\alpha)=(\lambda,\mu)$. Using Lemma \ref{thm-lem}, we find that $\lambda \in \mathcal{P}_{k-r}(k-s,m)$ and $\mu \in \mathcal{P}_{n-k+r}(0,\,k-s)$. We next show that $\lambda_1=m$. Recall that $\alpha^{(i)}$ denotes the permutation obtained by removing $\pi_1,\ldots,\pi_{i}$ from $\alpha$. From the definition  of $\Phi$ (that is, Definition \ref{defi-map}), we see that
\[\lambda_1=t(i_{k-r})={\rm maj}(\alpha^{(i_{k-r}-1)})-{\rm maj}(\alpha^{(i_{k-r})})-d_{i_{k-r}}(\pi),\]
where $i_{k-r}$ is the largest number such that ${\rm des}(\alpha^{(i_{k-r}-1)})={\rm des}(\alpha^{(i_{k-r})})+1$. Under the precondition that $\alpha \in \mathfrak{S}^{la}(\sigma,\pi|k)$, we know that there exists ${\rm ld}(\pi)\leq j\leq n$ (where ${\rm ld}(\pi)$ is the largest descent of $\pi$  if it exists and $0$ otherwise) satisfying the conditions (1)--(4) in Definition \ref{defi-refi-1}. It is easy to check that ${\rm des}(\alpha^{(j-1)})={\rm des}(\alpha^{(j)})+1$. Moreover, for $p>j$, we have ${\rm des}(\alpha^{(p-1)})={\rm des}(\alpha^{(p)})$. It implies that $i_{k-r}=j$.
Moreover,
\[{\rm maj}(\alpha^{(i_{k-r}-1)})-{\rm maj}(\alpha^{(i_{k-r})})-d_{i_{k-r}}(\pi)=m,\]
and so $\lambda_1=t(i_{k-r})=m$. Hence we derive that $\Phi(\alpha)=(\lambda,\mu) \in \mathcal{P}_{k-r}(k-s,=m)\times \mathcal{P}_{n-k+r}(0,\,k-s)$.

To prove that the  map $\Phi$ is the desired bijection, we proceed to show that the inverse map $\Psi$ given in Definition \ref{defi-psi} is  also  a map from
 $ \mathcal{P}_{k-r}(k-s, =m)\times \mathcal{P}_{n-k+r}(0,\,k-s)$ to $\mathfrak{S}^{la}(\sigma,\pi|k)$.  Let $\lambda \in \mathcal{P}_{k-r}(k-s,=m)$ and $\mu \in \mathcal{P}_{n-k+r}(0,\,k-s)$.
Assume that $\Psi(\lambda,\mu)=\overline{\alpha}$. We aim to show that $\overline{\alpha} \in \mathfrak{S}^{la}(\sigma,\pi)$.

With the help of Lemma \ref{lem2}, we know that $\overline{\alpha} \in \mathfrak{S}(\sigma,\pi|k)$. It remains to show that there exists $j$ satisfying the conditions in Definition \ref{defi-refi-1}.  Recall that $M^{(n)}$ is given by \eqref{part-cond-multiset}
  and
  $T^{(n-b)}$ is given by \eqref{pf-inverse-a} for $0\leq b\leq n$.  Let $k_{n-b}$ be the largest positive integer such that $T^{(n-b)}_{k_{n-b}}\in M^{(n-b)}$. Then  $M^{(n-b-1)}$ is obtained from $M^{(n-b)}$ by removing one $T^{(n-b)}_{k_{n-b}}$ and the permutation $\alpha^{(n-b-1)}$ is obtained   by inserting $\pi_{n-b}$ into $\alpha^{(n-b)}$ before  $\alpha_{k_{n-b}}^{(n-b)}$.

Since $\lambda_1=m$, we find that $m \in M^{(n)}$. From Corollary \ref{3.1}, we deduce that $m \in T^{(n)}$. Assume that $j$ is the largest number such that  $T_{k_j}^{(j)}=m$.  We proceed to prove that $j$ satisfies the conditions in the definition of $\mathfrak{S}^{la}(\sigma,\pi|k)$.

Since $j$  is the largest number such that  $T_{k_j}^{(j)}=m$, that means that for $n\geq p>j$, $T_{k_p}^{(p)}\leq k-s$, it implies that  ${\rm des}(\overline{\alpha}^{(p-1)})={\rm des}(\overline{\alpha}^{(p)})$, and so ${\rm des}(\overline{\alpha}^{(j)})={\rm des}(\overline{\alpha}^{(n)})={\rm des}(\sigma)$.
Observe that $\overline{\alpha}^{(j-1)}$ is obtained by inserting $\pi_j$ into  $\overline{\alpha}^{(j)}$ before  $\overline{\alpha}_{k_j}^{(j)}$, so all the elements in $\overline{\alpha}^{(j-1)}$ before $\pi_j$ are all from $\sigma$, that is, $\sigma_1,\ldots ,\sigma_{k_j-1}$. Hence it follows that
\begin{equation}\label{pf-r-tem}
d_{k_j-1}(\overline{\alpha}^{(j)})=d_{k_j-1}(\sigma)
\end{equation}
since ${\rm des}(\overline{\alpha}^{(j)})={\rm des}(\sigma)$.
   Recall that $T_{k_j}^{(j)}=m$, by definition,  we have
\begin{equation}\label{pf-m-temp}
{\rm maj}(\overline{\alpha}^{(j-1)})-{\rm maj}(\overline{\alpha}^{(j)})-d_{j}(\pi)=m.
\end{equation}
It means that ${\rm des}(\overline{\alpha}^{(j-1)})={\rm des}(\overline{\alpha}^{(j)})+1$, and the new generated descent of $\overline{\alpha}^{(j-1)}$ would be $k_{j}-1$ or $k_{j}$. Hence, from \eqref{pf-r-tem}, we deduce that  there are only three possible cases where \eqref{pf-m-temp} holds.
\begin{itemize}
\item[Case 1.]  $k_{j}-1$ is the new generated descent of $\overline{\alpha}^{(j-1)}$, $d_{k_j}(\overline{\alpha}^{(j)})=m-k_j+1$ and $d_{j}(\pi)=0$; or
\item[Case 2.]   $k_j$ is the new generated descent of $\overline{\alpha}^{(j-1)}$, $d_{k_j}(\overline{\alpha}^{(j)})=m-k_j$ and $d_{j}(\pi)=0$; or
\item[Case 3.]  $k_j$ is the new generated descent of $\overline{\alpha}^{(j-1)}$, $d_{k_j}(\overline{\alpha}^{(j)})=m-k_j+1$ and $d_{j}(\pi)=1$.
\end{itemize}
It can be checked that $j$ satisfies the conditions in the definition of $\mathfrak{S}^{la}(\sigma,\pi)$ under the condition that ${\rm des}(\overline{\alpha}^{(j)})={\rm des}(\overline{\alpha}^{(j+1)})=\cdots={\rm des}(\overline{\alpha}^{(n)})={\rm des}(\sigma)$.  More precisely, if $\overline{\alpha}^{(j-1)}$ satisfies the conditions stated in Case 1, it is easy to check that $j$ satisfies the conditions (1) and (2) in Definition \ref{defi-refi-1}. When $\overline{\alpha}^{(j-1)}$ belongs to Case 2, we find that $j$ satisfies
the conditions (1) and (4) in Definition \ref{defi-refi-1}. If $\overline{\alpha}^{(j-1)}$ meets the conditions in Case 3, then $j$ satisfies  the conditions (1) and (3) in Definition \ref{defi-refi-1}.

Since $\overline{\alpha}$ is obtained from $\overline{\alpha}^{(j-1)}$ by inserting $\pi_1,\ldots,\pi_{j-1}$ before $\pi_{j}$, it does not affect the aforementioned conclusion that $j$ satisfies the conditions in the definition of $\mathfrak{S}^{la}(\sigma,\pi)$. It follows that  $\Psi(\lambda,\mu)=\overline{\alpha}  \in \mathfrak{S}^{la}(\sigma,\pi)$. Hence we conclude that the map $\Psi$  is   a map from
 $ \mathcal{P}_{k-r}(k-s, =m)\times \mathcal{P}_{n-k+r}(0,\,k-s)$ to $\mathfrak{S}^{la}(\sigma,\pi|k)$.
Since the map $\Phi$ and the map $\Psi$ are inverse to each other, we deduce that the map $\Phi$ is the bijection as desired.
This completes the proof. \qed

Applying \eqref{int-GassCoeft}  and \eqref{int-GassCoeft-3}, as well as Theorem \ref{Stanley-c},   we find that the proof of Theorem \ref{stanley-r3} is equivalent to the proof of the following combinatorial assertion.

 \begin{thm}\label{Stanley-c-r3} Suppose that $\sigma \in \mathfrak{S}_m$ and $\pi \in \mathfrak{S}_n$ are two disjoint permutations, where ${\rm des}(\sigma)=r$  and  ${\rm des}(\pi)=s$. Then the  map $\Phi$ given in Definition \ref{defi-map} is   a bijection between $\mathfrak{S}^{sa}(\sigma,\pi|k)$ and $ \mathcal{P}_{k-r}(=k-s, m)\times \mathcal{P}_{n-k+r}(0,\,k-s)$.
 \end{thm}

\proof Let $\alpha \in \mathfrak{S}^{sa}(\sigma,\pi|k)$. Assume that $\Phi(\alpha)=(\lambda,\mu)$. According to Lemma \ref{thm-lem}, we see that $\lambda \in \mathcal{P}_{k-r}(k-s,m)$ and $\mu \in \mathcal{P}_{n-k+r}(0,\,k-s)$, namely,
\[  m\geq \lambda_1\geq \cdots \ge \lambda_{k-r}\geq k-s \geq  {\mu}_1\geq \cdots \geq {\mu}_{n-k+r}\geq 0. \]
We proceed to show that $\lambda_{k-r}=k-s$.
 Since $\alpha \in \mathfrak{S}^{sa}(\sigma,\pi|k)$,  there exists $1\leq j \leq {\rm sd}(\pi)$ (where ${\rm sd}(\pi)$ is the smallest descent of $\pi$ if it exists and $n$ otherwise) satisfying the conditions in Definition \ref{defi-refi-3}. Recall that $\alpha^{(i)}$ is the permutation obtained by removing $\pi_1,\ldots,\pi_i$ from $\alpha$. From Definition \ref{defi-refi-3}, we see that ${\rm des}(\alpha^{(j-1)})={\rm des}(\alpha^{(j)})+1$,  $d_j(\pi)={\rm des}(\pi)=s$, and for
 $p<j$, ${\rm des}(\alpha^{(p-1)})={\rm des}
 (\alpha^{(p)})$. Hence we conclude that ${\rm des}
 (\alpha^{(j-1)})={\rm des}(\alpha^{(0)})={\rm des}
 (\alpha)=k$, and by ${\rm des}(\alpha^{(j-1)})={\rm des}
 (\alpha^{(j)})+1$, we deduce that ${\rm des}
 (\alpha^{(j)})=k-1$. Since $j$ meets the conditions in Definition \ref{defi-refi-3}, we derive that elements before $\pi_{j}$ in  $\alpha^{(j-1)}$ are from $\sigma$, we may assume that $\sigma_1, \ldots, \sigma_q$. Then   $\sigma_1>\cdots>\sigma_q>\pi_{j}$, so we arrive at  ${\rm maj}(\alpha^{(j-1)})-{\rm maj}(\alpha^{(j)})=k$. Moreover, by the definition of $\Phi$, we see that
\[
\lambda_{k-r}=t(j)={\rm maj}(\alpha^{(j-1)})-{\rm maj}(\alpha^{(j)})-d_{j}(\pi)=k-s.
\]
It yields that  $\lambda \in \mathcal{P}_{k-r}(=k-s,m)$.

 We proceed to show that the inverse map $\Psi$ given in Definition \ref{defi-psi} is also  a map from
 $\mathcal{P}_{k-r}(=k-s,m)\times \mathcal{P}_{n-k+r}(0,\,k-s)$ to $\mathfrak{S}^{sa}(\sigma,\pi|k)$.    Let $\lambda \in \mathcal{P}_{k-r}(=k-s,m)$ and $\mu \in \mathcal{P}_{n-k+r}(0,\,k-s)$. Assume that $\Psi(\lambda,\mu)=\overline{\alpha}$.
From Lemma \ref{lem2}, we derive that $\overline{\alpha}=\overline{\alpha}_1\cdots \overline{\alpha}_{n+m}$  is a shuffle of $\sigma$ and $\pi$ with $k$ descents.  Assume that $\overline{\alpha}^{(i-1)}$ is the permutation obtained from $\overline{\alpha}^{(i)}$ by inserting $\pi_i$ before $\overline{\alpha}_{k_i+1}^{(i)}$ and  $j$ is the smallest integer such that ${\rm des}(\overline{\alpha}^{(j-1)})={\rm des}(\overline{\alpha}^{(j)})+1$.  We proceed to show that $j$ satisfies the conditions in  Definition \ref{defi-refi-3}.  Set $k_{j}=q$.
Then $\sigma_1,\ldots, \sigma_q$ must appear before $\pi_{j}$ in $\overline{\alpha}^{(j-1)}$. Under the assumption that $j$ is the smallest integer such that ${\rm des}(\overline{\alpha}^{(j-1)})={\rm des}(\overline{\alpha}^{(j)})+1$, we derive that ${\rm des}(\overline{\alpha}^{(j-1)})={\rm des}(\overline{\alpha})=k$, and so ${\rm des}(\overline{\alpha}^{(j)})=k-1$. Consequently,
${\rm maj}(\overline{\alpha}^{(j-1)})-{\rm maj}(\overline{\alpha}^{(j)})\geq k$. In particular, we deduce that if ${\rm maj}(\overline{\alpha}^{(j-1)})-{\rm maj}(\overline{\alpha}^{(j)})=k$, then $\sigma_1>\cdots>\sigma_q>\pi_{j}$ and either $\sigma_q<\overline{\alpha}_{q+1}^{(j)}$ or $\pi_{j}>\overline{\alpha}_{q+1}^{(j)}$.
From the proofs of Lemma \ref{thm-lem}
 and Lemma \ref{lem2}, we find that
 $\lambda_{k-r}=t(j)={\rm maj}(\overline{\alpha}^{(j-1)})-{\rm maj}(\overline{\alpha}^{(j)})-d_{j}(\pi).$
Since $\lambda_{k-r}=k-s$, we see that
\begin{equation}\label{pf-refine3}
{\rm maj}(\overline{\alpha}^{(j-1)})-{\rm maj}(\overline{\alpha}^{(j)})-d_{j}(\pi)=k-s.
\end{equation}
Observe that ${\rm maj}(\overline{\alpha}^{(j-1)})-{\rm maj}(\overline{\alpha}^{(j)})\geq k$ and $d_{j}(\pi)\leq s$, in order for \eqref{pf-refine3} to be valid, it's  necessary  that
\[{\rm maj}(\overline{\alpha}^{(j-1)})-{\rm maj}(\overline{\alpha}^{(j)})=k \quad \text{and} \quad  d_{j}(\pi)=s.\] Hence we conclude that  $j$ satisfies the condition in Definition \ref{defi-refi-3}, and so   $\overline{\alpha} \in \mathfrak{S}^{sa}(\sigma,\pi)$. This completes the proof.
\qed

We conclude this paper with the proof of Theorem \ref{stanley-r4}. In light of \eqref{int-GassCoeft}  and \eqref{int-GassCoeft-2}, together with  Theorem \ref{Stanley-c}, it is necessary to  prove the following combinatorial statement.

 \begin{thm}\label{Stanley-c-r4} Suppose that $\sigma \in \mathfrak{S}_m$ and $\pi \in \mathfrak{S}_n$ are two disjoint permutations, where ${\rm des}(\sigma)=r$  and  ${\rm des}(\pi)=s$. Then the  map $\Phi$ given in Definition \ref{defi-map} is   a bijection between $\mathfrak{S}^{lb}(\sigma,\pi|k)$ and $ \mathcal{P}_{k-r}(k-s, m)\times \mathcal{P}_{n-k+r}(0,=k-s)$.
 \end{thm}

 \pf  Let $\alpha \in \mathfrak{S}^{lb}(\sigma,\pi|k)$. Assume that $\Phi(\alpha)=(\lambda,\mu)$. From Lemma \ref{thm-lem}, we see that $\lambda \in \mathcal{P}_{k-r}(k-s,m)$ and $\mu \in \mathcal{P}_{n-k+r}(0,\,k-s)$, namely,
\[  m\geq \lambda_1\geq \cdots \ge \lambda_{k-r}\geq k-s \geq  {\mu}_1\geq \cdots \geq {\mu}_{n-k+r}\geq 0. \]
We proceed to show that ${\mu}_1=k-s$.
 Since $\alpha \in \mathfrak{S}^{lb}(\sigma,\pi|k)$,  there exists $1\leq j \leq {\rm sa}(\pi)$ (where ${\rm sa}(\pi)$ is the smallest ascent of $\pi$) satisfying the conditions in Definition \ref{defi-refi-4}.   Recall that $\alpha^{(i)}$ is the permutation obtained by removing $\pi_1,\ldots,\pi_i$ from $\alpha$. From Definition \ref{defi-refi-4}, we see that ${\rm des}(\alpha^{(j-1)})={\rm des}(\alpha^{(j)})$,  $d_{j}(\pi)=s-j+1$. Moreover, the elements before $\pi_j$ in $\alpha^{(j-1)}$ are from $\sigma$, that is, $\sigma_1<\cdots<\sigma_q$. And for
 $a<j$, ${\rm des}(\alpha^{(a-1)})={\rm des}
 (\alpha^{(a)})+1$. Hence we conclude that ${\rm des}
 (\alpha^{(j-1)})={\rm des}(\alpha^{(0)})-j+1=k-j+1$, and by ${\rm des}(\alpha^{(j-1)})={\rm des}
 (\alpha^{(j)})$, we deduce that ${\rm des}
 (\alpha^{(j)})=k-j+1$. Under the condition that the elements before $\pi_{j}$ in  $\alpha^{(j-1)}$ are from $\sigma$, that is, $\sigma_1<\cdots<\sigma_q$, we deduce that   ${\rm maj}(\alpha^{(j-1)})-{\rm maj}(\alpha^{(j)})=k-j+1$. Moreover, by the definition of $\Phi$, we see that
\[
\mu_{1}=t(j)={\rm maj}(\alpha^{(j-1)})-{\rm maj}(\alpha^{(j)})-d_{j}(\pi)=k-j+1-(s-j+1)=k-s.
\]
It yields that  $\mu \in \mathcal{P}_{n-k+r}(0,=k-s)$.

 Conversely,   let $\lambda \in \mathcal{P}_{k-r}(k-s,m)$ and $\mu \in \mathcal{P}_{n-k+r}(0,=k-s)$, we next show that the map $\Psi$ given in Definition \ref{defi-psi} is also  a map from
 $\mathcal{P}_{k-r}(k-s,m)\times \mathcal{P}_{n-k+r}(0,=k-s)$ to $\mathfrak{S}^{lb}(\sigma,\pi|k)$.  Assume that $\Psi(\lambda,\mu)=\overline{\alpha}$. First, we can use Lemma \ref{lem2} to derive that $\overline{\alpha}=\overline{\alpha}_1\cdots \overline{\alpha}_{n+m}$  is a shuffle of $\sigma$ and $\pi$ with $k$ descents.  Assume that $\overline{\alpha}^{(i-1)}$ is the permutation obtained from $\overline{\alpha}^{(i)}$ by inserting $\pi_i$ before $\overline{\alpha}_{k_i+1}^{(i)}$ and  $j$ is the smallest integer such that ${\rm des}(\overline{\alpha}^{(j-1)})={\rm des}(\overline{\alpha}^{(j)})$.  We aim to show that $j$ satisfies the conditions in  Definition \ref{defi-refi-4}.  Set $k_{j}=q$, it means that $\sigma_1,\ldots, \sigma_q$ should appear before $\pi_{j}$ in $\overline{\alpha}^{(j-1)}$. Under the assumption that $j$ is the smallest integer such that ${\rm des}(\overline{\alpha}^{(j-1)})={\rm des}(\overline{\alpha}^{(j)})$, we derive that ${\rm des}(\overline{\alpha}^{(j-1)})=k-j+1$. Consequently,
${\rm maj}(\overline{\alpha}^{(j-1)})-{\rm maj}(\overline{\alpha}^{(j)})\leq k-j+1$. In particular, we deduce that if ${\rm maj}(\overline{\alpha}^{(j-1)})-{\rm maj}(\overline{\alpha}^{(j)})=k-j+1$, then $\sigma_1<\cdots<\sigma_q<\pi_{j}$ and either $\sigma_q>\overline{\alpha}^{(j)}_{q+1}$ or $\pi_{j}<\overline{\alpha}^{(j)}_{q+1}$.
From the proofs of Lemma \ref{thm-lem}
 and Lemma \ref{lem2}, we find that $\mu_{1}=t(j)={\rm maj}(\overline{\alpha}^{(j-1)})-{\rm maj}(\overline{\alpha}^{(j)})-d_{j}(\pi).$
Since $\mu_{1}=k-s$, we see that
\begin{equation}\label{refin4-pfa}
{\rm maj}(\overline{\alpha}^{(j-1)})-{\rm maj}(\overline{\alpha}^{(j)})-d_{j}(\pi)=k-s.
\end{equation}
Since ${\rm maj}(\overline{\alpha}^{(j-1)})-{\rm maj}(\overline{\alpha}^{(j)})\leq k-j+1$ and $d_{j}(\pi)\geq s-j+1$, we find that \eqref{refin4-pfa} holds if and only if
\[{\rm maj}(\overline{\alpha}^{(j-1)})-{\rm maj}(\overline{\alpha}^{(j)})=k-j+1 \quad \text{and}  \quad d_{j}(\pi)=s-j+1.\]
Hence we conclude that  $j$ satisfies the condition in Definition \ref{defi-refi-4}, and so   $\overline{\alpha} \in \mathfrak{S}^{lb}(\sigma,\pi)$. This completes the proof.
\qed

 \vskip 0.2cm
\noindent{\bf Acknowledgment.} This work
was supported by   the National Science Foundation of China.


\begin{thebibliography}{0}
	\setlength{\itemsep}{-.8mm}
	\addcontentsline{toc}{section}{}


  \bibitem{Adin-Gessel-Reiner-Roichman-2021} R.M. Adin, I.M. Gessel, V. Reiner, and Y. Roichman, Cyclic quasi-symmetric functions, Israel J. Math. 243 (2021) 437--500.
	
 \bibitem{Andrews-1976} G.E. Andrews, The Theory of Partitions,
 Addison-Wesley Publishing Co., 1976.

 \bibitem{Baker-Jarvis-Sagan-2020} D. Baker-Jarvis and B.E. Sagan, Bijective proofs of shuffle compatibility results, Adv. in Appl. Math. 113 (2020) 101973.


 \bibitem{Domagalski-Liang-Minnich-Sagan-2020} R. Domagalski, J. Liang, Q. Minnich, B.E. Sagan, J. Schmidt and A. Sietsema, Cyclic shuffle compatibility, S\'em. Lothar. Combin. 85 ([2020--2021]), Art. B85d, 11 pp.

 \bibitem{Foata-1968} D. Foata, On the Netto inversion number of a sequence, Proc. Amer. Math. Soc. 19 (1968) 236--240.

 \bibitem{Garsia-Gessel-1979} A.M. Garsia and I.M. Gessel, Permutation statistics and partitions,  Adv.   Math. 31 (1979) 288--305.

  \bibitem{Gessel-Zhuang-2018} I.M. Gessel and Y. Zhuang, Shuffle-compatible permutation statistics, Adv. Math. 332 (2018) 85--141.

   \bibitem{Grinberg-2018} D. Grinberg, Shuffle-compatible permutation statistics II: the exterior peak set,
Electron. J. Combin. 25 (2018) P4.17.

\bibitem{Goulden-1985}  I. P. Goulden, A bijective proof of Stanley's shuffling theorem, Trans. Amer. Math. Soc. 288 (1985)  147--160.

\bibitem{Haglund-Loehr-Remmel-2005} J. Haglund, N. Loehr and J.B. Remmel, Statistics on wreath products, perfect matchings, and signed words,
European J. Combin. 26 (2005)  835--868.

\bibitem{Ji-Zhang-2022} K.Q. Ji and D.T.X. Zhang, A cyclic analogue of Stanley's shuffling theorem, Electron. J. Combin. 29 (2022) P4.20.

\bibitem{MacMahon-1978} P.A. MacMahon, Collected Papers. Vol. I, Cambridge, MA, 1978.

\bibitem{Novick-2010} M. Novick, A bijective proof of a   major index theorem of Garsia and Gessel, Electron. J. Combin. 17 (2010), Research Paper 64, 12 pp.

\bibitem{Sagan-Savage-2012} B.E. Sagan and C.D. Savage, Mahonian pairs, J. Combin. Theory Ser. A 119 (2012) 526--545.

\bibitem{Stadler-1999} J. D. Stadler, Stanley's shuffling theorem revisited, J. Combin. Theory Ser. A 88 (1999)  176--187.




\bibitem{Stanley-1972} R. P. Stanley, Ordered structures and partitions, Mem Amer. Math. Soc., No. 119. American Mathematical Society, 1972. iii+104 pp.

 \bibitem{Stanley-2012}   R.P. Stanley, Enumerative Combinatorics, Vol. I, 2nd ed.,  Cambridge University Press, Cambridge, 2012.

\bibitem{Yang-Yan-2022} L. Yang and S.H.F. Yan, On a conjecture concerning shuffle-compatible permutation statistics, Electron. J. Combin. 29 (2022) P3.3.


\end{thebibliography}
\end{document}